\documentclass[12pt,reqno]{amsart}
\usepackage{amssymb}
\usepackage{geometry}
\usepackage[latin1]{inputenc}
\usepackage[italian,english]{babel}
\usepackage{amsmath, amsfonts, amsthm}
\usepackage{graphicx}
\usepackage{pgfplots}
\usepackage{paralist}
\usepackage{mathtools, mathabx,bigints}
\usepackage{subfig}
\usepackage{comment}
\usepackage{adjustbox}
\usepackage{cleveref}
\usepackage{siunitx}
\usepackage{soul}
\usepackage{algpseudocodex}
\usepackage{tcolorbox} 
\numberwithin{equation}{section}
\newtheorem*{definition}{Definition}

\newcommand{\E}{\mathbf{E}}
\newcommand{\DD}{\mathbf{D}}

\newcommand{\II}{\mathbf{I}}

\newcommand{\J}{\mathbf{J}}
\newcommand{\mcR}{\mathbb{R}}

\newcommand{\QQ}{\mathbf{Q}}
\newcommand{\rr}{\mathbf{r}}
\newcommand{\RR}{\mathbf{R}}
\newcommand{\w}{\mathbf{w}}
\newcommand{\xx}{\mathbf{x}}

\newcommand{\ZZ}{\mathbf{Z}}

\newcommand{\HD}{H(\text{div};\Omega)}

\let\oldoint\oint
\renewcommand{\oint}{\oldoint\nolimits}

{\left\{\begin{array}{@{}l@{}}}{\end{array}\right.}
\patchcmd{\abstract}{\scshape\abstractname}{\textbf{\abstractname}}{}{}
\makeatletter 
\def\@makefnmark{} 
\makeatother 
 
\pgfplotsset{compat=1.18}
\begin{document}
\title[QR COMPRESSION OF VIE FOR METASURFACES]{QR-Recursive Compression of Volume Integral Equations for Electromagnetic Scattering by Large Metasurfaces}
\author[\,]{
Vincenzo Mottola$^1$, Antonello Tamburrino$^1$, Luca Bergamaschi$^2$, Andrea G. Chiariello$^3$, Emanuele Corsaro$^4$, Carlo Forestiere$^4$, Guglielmo Rubinacci$^{4,5}$ and Salvatore Ventre$^1$}\footnote{\\$^1$Department of Electrical and Information Engineering, University of Cassino and Southern Lazio, Via Gaetano di Biasio 43, 03043 Cassino, Italy.\\
$^2$ Department of Civil Environmental and Architectural Engineering University of Padova,  via Marzolo 9, 35131 Padova \\
$^3$ Department of Engineering, Università degli Studi della Campania \lq\lq Luigi Vanvitelli\rq\rq, via Roma 29, 81031 Aversa, Italy \\
$^4$ Department of Electrical Engineering and Information Technology,
Università degli Studi di Napoli Federico II, 80138 Napoli, Italy \\
$^5$ Consorzio CREATE, 80125 Napoli, Italy \\
Email: antonello.tamburrino@unicas.it }

\begin{abstract} 
In this paper, a novel QR decomposition-based compression scheme is combined with a  volume integral equations method for the fast and efficient numerical computation of the scattering of electromagnetic fields from large scale metasurfaces, via an iterative approach.

The underlying problem is of a multiscale nature. Indeed, these metasurfaces are made of a large collection of interacting sub-wavelength scatterers, thus making the numerical computation of the solution very challenging.

More specifically, the paper proposes a tailored version of a QR decomposition-based compression for a volume integral equation, together with a proper preconditioner that exploits the geometrical structure of the array, in order to achieve a fast and accurate iterative solver, in view of realistic applications.  Numerical examples prove the effectiveness of the method in efficiently modeling metasurfaces made by thousands of particles.

\end{abstract}

\maketitle
\markright{QR COMPRESSION OF VIE FOR METASURFACES}

\section{Introduction} \label{sec:introduction}
Metasurfaces, composed of interacting subwavelength scatterers operating near resonance \cite{holloway_overview_2012}, offer a powerful platform for controlling electromagnetic waves. By carefully designing the geometry and arrangement of these scatterers, metasurfaces enable the engineering of a wide range of effective constitutive relations \cite{pfeiffer_metamaterial_2013}. Their versatility has driven advances in applications such as analog computing \cite{silva_performing_2014}, biosensing \cite{altug_advances_2022}, and imaging \cite{chen_review_2016}. Efficient and accurate simulation of metasurfaces is essential for understanding their fundamental properties \cite{choi_realization_2024} and for optimizing their performance in practical implementations \cite{kang_large-scale_2024}.

Metasurfaces design often relies on the \lq\lq unit cell\rq\rq\ approximation \cite{yu_flat_2014,aieta_multiwavelength_2015,lalanne_prehistory_2023}. This approach constructs a library of the electromagnetic response of meta-atoms, such as disks, pillars, or nanofins, as a function of a design parameter, like the disk radius, pillar height, or the nanofin orientation, assuming the structure to be locally periodic. Then, the contributions of the different \lq\lq unit cell\rq\rq\ are combined by linear superposition. Unfortunately, the \lq\lq unit cell\rq\rq\ approach has several limitations: the most critical limitation is the assumption that the interaction between adjacent meta-atoms is assumed equal to that in periodic arrays, which necessitates a slowly changing structure. This constraint limits the design space, resulting in narrow-bandwidth operation, low efficiencies, and many types of non-ideality.

General-purpose electromagnetic simulation tools, based on either \textit{integral} or \textit{differential} formulations, face significant challenges in the simulation of metasurfaces \cite{choi_realization_2024}. Differential formulations, including finite-element methods and finite-difference methods, are easy to implement and result in sparse matrices but require truncation of the computational domain through the introduction of convenient boundary conditions \cite{jin_theory_2011}.   Among them, the finite-difference time-domain (FDTD) method \cite{yee_finite-difference_1997,taflove_computational_2005,oskooi_meep_2010}  is widely used to model metasurfaces. However, the finite resolution leads to inaccuracies in wave speed, resulting in phase accumulation errors as the size of the scattering region increases. Hence, even though the complexity of FDTD simulations nominally increases linearly with the size of the array, maintaining a fixed error requires increasing the resolution \cite{shlager_selective_1995,taflove_computational_2005}. This makes the approach prohibitive when applied using traditional computing platforms, requiring massive GPU parallelization to be viable \cite{10.1063/5.0071245}.

Volume Integral formulations (e.g. \cite{Rubinacci20062977, Miano20102920}) are appealing because they define the unknowns only within the objects' volumes, and naturally satisfy radiation conditions at infinity. In addition, contrarily to \textit{surface} formulations \cite{forestiere_static_2023,corsaro_multilevel_2024} they also allow to treat inhomogenous meta-atoms. However, they usually involve large dense matrices, which are computationally demanding to store and invert. Moreover, these matrices are often ill-conditioned, which can hinder the efficiency of iterative solvers such as GMRES.

To address these limitations, in this work we combine the volume integral formulation presented in \cite{Rubinacci20062977, Miano20102920} with a QR compression scheme and a dedicated preconditioner. QR factorization methods are not new in computational electromagnetics, having been successfully applied to surface integral equation methods \cite{poirier_numerically_1998} and volume integral magnetoquasistatic formulations \cite{RUBINACCI2008,art:qrspars}.
In this case, both the compression scheme and the preconditioner are specifically designed for multiscale collections of interacting meta-atoms, grouping basis functions by meta-atom at the finest-level group.
The proposed method represents a first and key contribution in view of the simulation of suitable large-scale electromagnetic problems where the scatterer extends for thousands of wavelengths. Indeed, this method can be easily coupled with Butterfly-Acceleration methods to tackle large-scale electromagnetic problems \cite{Bin2022,Guo2013}.

The paper is organized as follows: in \Cref{matmod} the mathematical formulation of VIE for the problem of interest is presented; in \Cref{nummod} the numerical counterpart of the continuous model is detailed; in \Cref{sec_prec} the iterative method and the pre-conditioner employed in the computations are presented, while in \Cref{sec_qr} the compression method is detailed. Finally, in \Cref{examples} some numerical examples on realistic configurations are reported and in \Cref{conclusions} the conclusions are drawn.

\section{Mathematical model} \label{matmod}
In this Section the mathematical model represented by the source integral volume formulation is briefly summarized (see \cite{Rubinacci20062977,Miano20102920,Forestiere20171224} for further details).

The underlying integral formulation is expressed in terms of equivalent sources. Specifically, the unknown is the polarization current density $\J$ induced in the dielectrics. It is assumed that: (i) the region $\Omega$ occupied by the dielectrics consists of one or more bounded Lipschitz domains of $\mcR^3$, (ii) dielectric materials are linear and non-spatially dispersive, (iii) the dielectric materials are homogeneous in each connected component of $\Omega$, (iv) magnetic materials are not present, (v) the fields are computed in time-harmonic operations at the angular frequency $\omega$ (the $e^{j \omega t}$ time behavior is assumed) and (vi) no source charges are present in any possible cavity of $\Omega$. 

In addition, for the sake of generality, conductive materials may be present. In this latter case, $\Omega$ is the region occupied by the conducting materials, too, the conducting materials are linear, non-spatially dispersive, and homogeneous in each connected component of $\Omega$.
In typical dielectric materials the electrical conductivity is vanishing, while the dielectric permittivity may have a non-vanishing imaginary part, in case of dielectric losses.

Hereafter, $\Sigma_i$ represents the $i-$th connected component of $\partial \Omega$.

The unknown of the integral formulation is, therefore, the sum of the polarization current density and of the (electrical) conduction current density, i.e.
\begin{equation}
\label{eq:const_rel}
    \J(\rr) = \left[\sigma(\rr)+j\omega \varepsilon_{0} \chi(\rr)\right] \E (\rr), \text{ in } \Omega,
\end{equation}
where $\sigma$ is the electrical conductivity, $\chi$ is the dielectric susceptibility, $\varepsilon_0$ is the dielectric permittivity of the free-space
 
In general, $\J$ is an element of $\HD$, the functional space of div-conforming vector fields:
\begin{equation}
    \HD=\{\mathbf{v} \in \mathbb{L}^2\left(\Omega \right): \nabla \cdot \mathbf{v} \in L^2\left(\Omega \right)\}.
\end{equation}
In the specific case, $\J$ is also a divergence free vector and its net flux on each connected component $\Sigma_i$ of the boundary is vanishing, because materials are homogeneous in each connected component of $\Omega$ and the cavities of $\Omega$ are free from source charges (see assumptions (iii) and (vi)). Thus, $\J$ is an element of the following solution space:
\begin{equation}
    H=\{\mathbf{v} \in \HD: \nabla \cdot \mathbf{v} = 0 \text{ in } \Omega, \int_{\Sigma_i} \mathbf{v}
\cdot\mathbf{\hat{n}} \text{d}S = 0, \, \forall i\}.
\end{equation}

The unknown $\J$ satisfies the following integral equation:
\begin{align}
 \left[\sigma(\rr)+j\omega\chi(\rr)\varepsilon_{0}\right]^{-1}\mathbf{J}\left(  \mathbf{r} \right)  +j\omega\mu_{0}\int_{\Omega}\mathbf{J}\left(  \mathbf{r}^{\prime}\right)  g\left(  \mathbf{r-r}^{\prime}\right)  \text{d} \tau ^{\prime} \nonumber  \\ +\frac{1}{j\omega\varepsilon_{0}}\nabla\int_{\partial \Omega}\left(  \mathbf{J}
\cdot\mathbf{\hat{n}}\right)  \left(  \mathbf{r}^{\prime}\right)  g\left(
\mathbf{r-r}^{\prime}\right)  \text{d}S^{\prime} =\mathbf{E}_{0}\left(  \mathbf{r}\right)
,\text{ in }\Omega
\label{eq:IntEq}
\end{align}
In equation \eqref{eq:IntEq} $\mu_0$ is the magnetic permeability of the free-space, $g\left(  \mathbf{r}\right)
\triangleq e^{-jkr}/\left(  4\pi r\right)$ is the scalar free space Green's function, $k=\omega \sqrt{\varepsilon_0 \mu_0}$ is the free-space wavelength, and $\E_0$ is the applied electric field.

Equation \eqref{eq:IntEq} in the weak form is:
\begin{align}
\label{eqnWeak}
& \int_{\Omega}  \J^{\prime
}(\rr) \cdot \left[\sigma(\rr)+j\omega\chi(\rr)\varepsilon_{0}\right]^{-1}\mathbf{J}\left(  \mathbf{r} \right) \text{d} \tau \nonumber\\
+ & j\omega\mu_{0} \int_{\Omega} \int_{\Omega} \J'(\rr) \cdot \mathbf{J}\left(  \mathbf{r}^{\prime
}\right)  g\left(  \mathbf{r-r}^{\prime}\right)  \text{d} \tau ^{\prime} \text{d} \tau \nonumber \\
+ & \frac{1}{j\omega\varepsilon_{0}} \int_{\partial \Omega} \int_{\partial \Omega} \left( \J^{\prime} \cdot \mathbf{\hat{n}} \right) (\rr) \left(  \mathbf{J}
\cdot \mathbf{\hat{n}}\right)  \left(  \mathbf{r}^{\prime}\right)  g\left(
\mathbf{r-r}^{\prime}\right)  \text{d}S^{\prime} \text{d}S \nonumber \\
= & \int_{\Omega} \J'(\rr) \cdot \mathbf{E}_{0}\left(  \mathbf{r}\right) \text{d} \tau,
\end{align}
where $\J^{\prime} \in H$ is an arbitrary test function.

It is worth noting that \eqref{eq:IntEq} and \eqref{eqnWeak} are also valid in case of anisotropic materials, i.e. when $\sigma+j\omega\chi\varepsilon_{0}$ is a tensor field.

\section{Numerical model} \label{nummod}
The numerical model is obtained by means of the Galerkin method, once the unknown vector field $\J$ has been properly discretized. Specifically, $\J$ is discretized in term of linear combination of loop and star shape functions, as described in \cite{Rubinacci20062977,Miano20102920,Forestiere20171224,Miano2007586,DalNegro20091618}:
\begin{equation}
\mathbf{J=}\sum_{k=1}^{N_{L}}I_{k}^{L}\mathbf{w}_{k}^{L}+\sum_{k=1}^{N_{S}%
}I_{k}^{S}\mathbf{w}_{k}^{S}.
\label{eqn04}%
\end{equation}
where the $\w_k^L$s and the $\w_k^S$s are the \emph{loop} and \emph{star} shape functions, respectively, $N_{L}$ and $N_{S}$ are the number of loop and star shape functions and $I_k^L$ and $I_k^S$ are the corresponding Degrees of Freedom (DoFs).

\begin{figure}[!ht]
    \centering
    \includegraphics[width=0.75\linewidth]{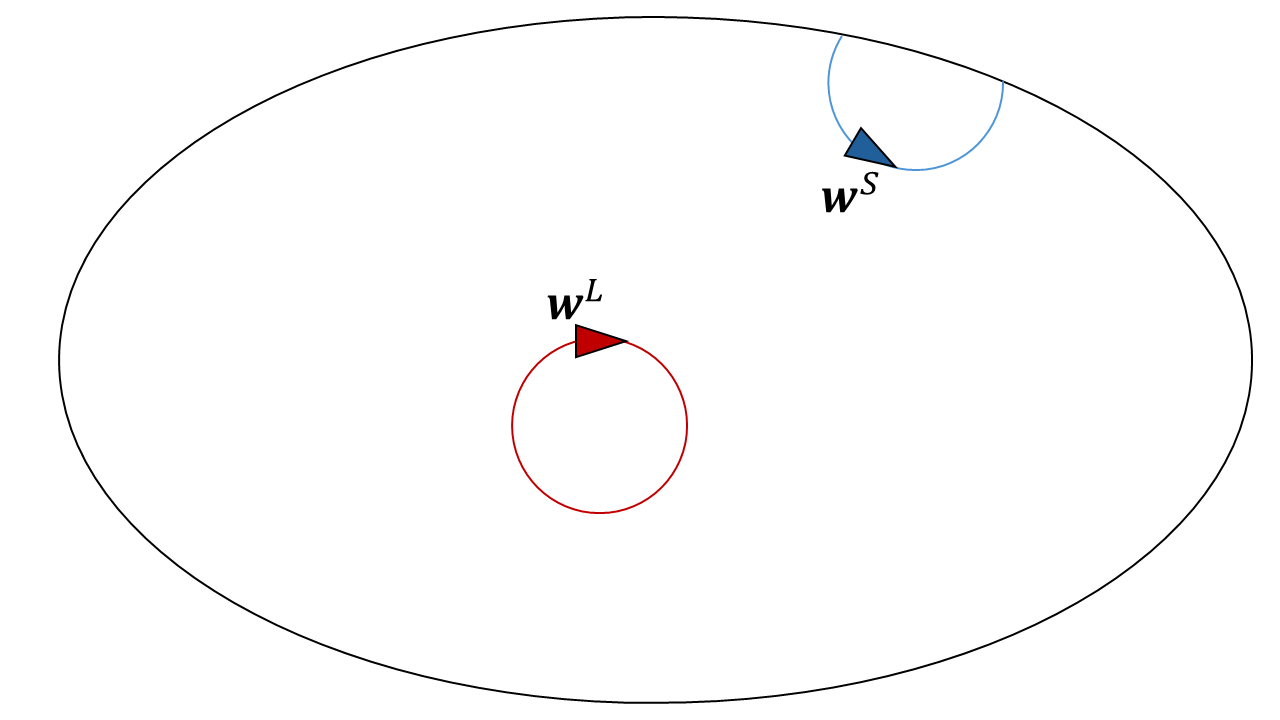}
    \caption{Examples of loop (orange) and star (blue) basis functions for a single meta-atom.}
    \label{fig01}
\end{figure}

When $\Omega$ is simply connected, a loop shape function is related to the edge of the finite element mesh and is given by the curl of the corresponding edge shape function \cite{Rubinacci20062977}. Each loop shape function represents an \emph{elementary} current loop. Star shape functions are related to the edges on boundary of the finite element mesh and are also given by the curl of the corresponding edge element. Star shape function corresponds to local \emph{truncated} loops (see Figure \ref{fig01}). Uniqueness of the representation \eqref{eqn04} can be guaranteed through graph theoretical approaches \cite{Rubinacci20062977,Forestiere20171224}.

The loop and star shape functions satisfy the following properties:
\begin{align}
\nabla \cdot \w_k^L  & =0\text{ in }\Omega\text{, }\w_k^L%
\cdot\mathbf{\hat{n}}=0\text{ on }\partial\Omega\\
\nabla \cdot \w_k^S  & =0\text{ in }\Omega, \ \w_k^S \cdot \mathbf{\hat{n}} \neq 0\text{ on }\partial\Omega, \ \int_{\Sigma_i} \w_k^S
\cdot\mathbf{\hat{n}} \text{d}S = 0, \, \forall i
\end{align}

The decomposition of the unknown in terms of loop and star components is mandatory to mitigate the so-called low frequency breakdown problem, that is, the dramatic increase of the condition number of the stiffness matrix in the low frequency limit and in the presence of conducting materials that prevent getting an accurate solution \cite{Rubinacci20062977,Zhao20001635}.

By plugging approximation \eqref{eqn04} in 
\eqref{eqnWeak} and applying the Galerkin method (projection along the $\w_k^L$s and the $\w_k^S$), it results that:
\begin{equation}
\label{matrixeq}
\left[ \mathbf{R} +j \omega \mu _{0}\mathbf{L} + \frac{1}{j\omega \varepsilon _{0}}\mathbf{D}\right] \mathbf{I}=\mathbf{v}_{0}
\end{equation}
where%
\begin{equation*}
\mathbf{R} = \left[ 
\begin{array}{cc}
\mathbf{R}^{LL} & \mathbf{R}^{LS} \\ 
\mathbf{R}^{SL} & \mathbf{R}^{SS}
\end{array}
\right] ;\ 
\mathbf{L} = \left[ 
\begin{array}{cc}
\mathbf{L}^{LL} & \mathbf{L}^{LS} \\ 
\mathbf{L}^{SL} & \mathbf{L}^{SS}
\end{array}
\right] ,
\end{equation*}
\begin{equation*}
\mathbf{D} =\left[ 
\begin{array}{cc}
\mathbf{0} & \mathbf{0} \\ 
\mathbf{0} & \mathbf{D}^{SS}
\end{array}
\right] ;\
\mathbf{v}_{0} = \left[ 
\begin{array}{c}
\mathbf{v}_{0}^{L} \\ 
\mathbf{v}_{0}^{S}
\end{array}
\right],
\end{equation*}
and
\begin{eqnarray*}
R_{ij}^{\alpha \beta } &=&\int_{\Omega }\mathbf{w}_{i}^{\alpha }\cdot \left[\sigma(\rr)+j\omega\chi(\rr)\varepsilon_{0}\right]^{-1} \mathbf{w}_{j}^{\beta }\text{d}V \\
L_{ij}^{\alpha \beta } &=&
\int_{\Omega} \int_{\Omega} \mathbf{w}_{i}^{\alpha }(\rr) \cdot g\left(  \mathbf{r-r}^{\prime}\right)  \mathbf{w}_{j}^{\beta }\left(  \mathbf{r}^{\prime
}\right) \text{d}V^{\prime} \text{d}V \\
v_{0,i}^{\alpha } &=&\int_{\Omega }\mathbf{w}_{i}^{\alpha }\cdot \mathbf{E}%
_{0}\text{d}V \\
D_{ij}^{SS} & = & 
\int_{\partial \Omega} \int_{\partial \Omega} \left( \mathbf{w}_{i}^{S} \cdot \mathbf{\hat{n}} \right) (\rr) g\left(
\mathbf{r-r}^{\prime}\right) \times \\
& & \times \left(\mathbf{w}_{j}^{S}
\cdot \mathbf{\hat{n}}\right)  \left(  \mathbf{r}^{\prime}\right) \text{d}S^{\prime} \text{d}S.
\end{eqnarray*}

In the following, we define
\begin{equation}
    \mathbf{Z}=\mathbf{R} +j \omega \mu _{0}\mathbf{L} + \frac{1}{j\omega \varepsilon _{0}}\mathbf{D}. 
\end{equation}
and, hence, equation \eqref{matrixeq} reads as 
\begin{equation}
\label{eq_vie}
    \mathbf{Z\,I}=\mathbf{v}_0.
\end{equation}
It is worth noticing that the $\mathbf{Z}$ matrix is fully populated.

Hereafter, the volume integral equation \eqref{eq_vie} is solved using an iterative method combining a QR compression method and a preconditioner specifically tailored for metasurfaces.

\section{Iterative method and preconditioner} \label{sec_prec}
When dealing with large-scale problems, direct solvers are impractical, leaving room for iterative solvers, representing the only way to solve the problem. Iterative solvers aim to compute the solution of \eqref{eq_vie} via the repeated computation of the $\mathbf{Z} \, \mathbf{I}_n$ and $\mathbf{Z}^H \, \mathbf{I}_n$ products, until the residual falls below a specified threshold.

\begin{figure}[h!]
    \centering
    \includegraphics[width=0.85\linewidth]{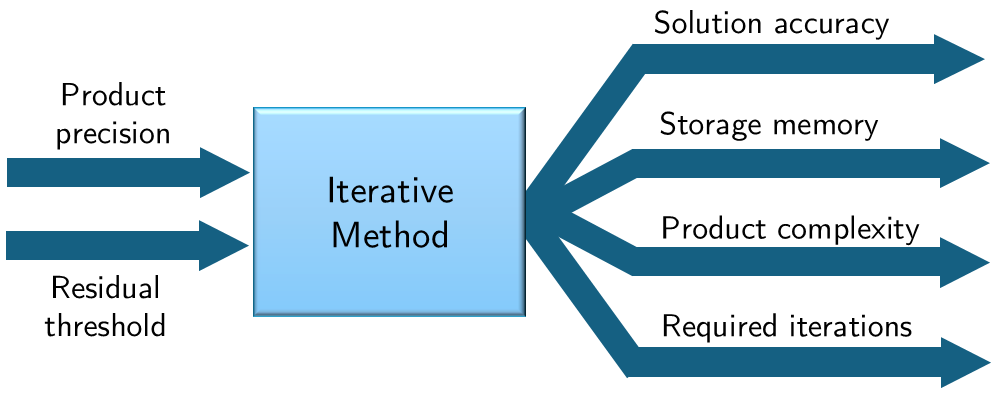}
    \caption{Iterative solvers.}
    \label{fig:qrbb}
\end{figure}

For large-scale problems, it is unpractical to store the full matrix $\mathbf{Z}$, so an approximate version obtained by a proper compression method is employed.
An iterative method can be conceptualized as a black box (Figure \ref{fig:qrbb}) that approximates the solution of a linear system. Given a specified precision for the compressed product $\mathbf{Z\,I}$ and a residual threshold, several key outputs are obtained. The first is solution accuracy, which quantifies the error between the numerically computed solution and the ``true'' solution, that could be determined using a direct method. Another important factor is storage memory, referring to the amount of memory required to evaluate the approximate $\mathbf{Z\,I}$ product at the given precision. Additionally, the method's efficiency is characterized by product complexity, which represents the computational cost of evaluating the $\mathbf{Z\,I}$ product. Finally, the number of iterations required to reduce the residual below the threshold provides insight into the convergence behavior of the method.

In reducing the computational cost, the preconditioners play a key role, since they aim to achieve the target residual value in a smaller number of iterations.

In this work, it is proposed an efficient and effective preconditioner tailored for the geometrical structure of the meta-atoms array. Specifically, the preconditioner consists of the exact inverse of the interaction matrix $\mathbf{Z}$ where the mutual interactions between different meta-atoms are neglected. 
For the sake of clarity, let the matrix $\ZZ$ be rearranged and partitioned by meta-atoms, i.e.
\begin{equation}
\label{bp}
\ZZ = \left[ 
\begin{array}{ccccc}
\ZZ_{1,1} & \ZZ_{1,2} & \ldots &  &\ZZ_{1,N} \\ 
\vdots & & & & \vdots \\
\ZZ_{N,1} & \ZZ_{N,2} & \ldots &  &\ZZ_{N,N} \\ 
\end{array}
\right],  
\end{equation}
where the $\ZZ_{ij}$ refers to the mutual interaction between meta-atom $i$ and $j$, and $N$ is the total number of meta-atoms. The block partition of \eqref{bp} is obtained by ordering the unknowns by meta-atoms, i.e. first the unknowns for the first meta-atom, then the unknowns for the second meta-atom, and so forth.

In this way, the (left) preconditioner used for solving \eqref{eq_vie} iteratively is the inverse of the block diagonal of $\ZZ$, i.e. is the inverse of
\begin{equation}
\label{diagprecond}
    \ZZ_D = diag\left(\ZZ_{1,1}, \, \ZZ_{2,2}, \, \ldots, \, \ZZ_{N,N} \right).
\end{equation}

This preconditioner solves exactly \eqref{eq_vie} when the interactions between meta-atoms are negligible, i.e. the meta-atoms are isolated.

The computational cost for evaluating this preconditioner increases only linearly with the total number of meta-atoms $N$, and reduces to $O(1)$ if all the meta-atoms are identical. Indeed, in this latter case, it is possible to use the same finite element mesh for each meta-atoms and, therefore, it results that $\ZZ_{1,1}= \ldots =\ZZ_{N,N}$. Moreover, the computational cost for computing the preconditioner depends on the number $N_p$ of DOFs present in one particle, only, and it results that 
(i) the number of operations for the factorization of the preconditioner is  $O(N_p^3)$,
(ii) the memory required to store the factorization of the preconditioner is $O(N_p^2)$, and
(iii) the number of operations to evaluate the preconditioner in the iterative scheme is $O(N_p^2)$.

\section{The QR-recursive method for a meta-atoms array }\label{sec_qr}

The core of the iterative method lies in approximating the $\ZZ \, \II$ matrix-by-vector product in a fast and accurate manner.

The proposed approach is based on the QR-recursive method, which has been successfully applied to MQS problems \cite{RUBINACCI2008, art:qrspars}, and here is adapted to a metasurface made by possibly many meta-atoms.
In the following, we briefly summarize the main aspects of the method, highlighting its novelty when applied to an array of meta-atoms.

The key concepts of the method are (i) the distinction of near and far distance interaction, (ii) the efficient compression via a low-rank QR-decomposition of the far distance interactions, and (iii) a recursive subdivision scheme.

\subsection{Far and near distance interactions}
\begin{figure}[h!]
    \centering
    \subfloat[][]
    {\includegraphics[width=0.45\linewidth]{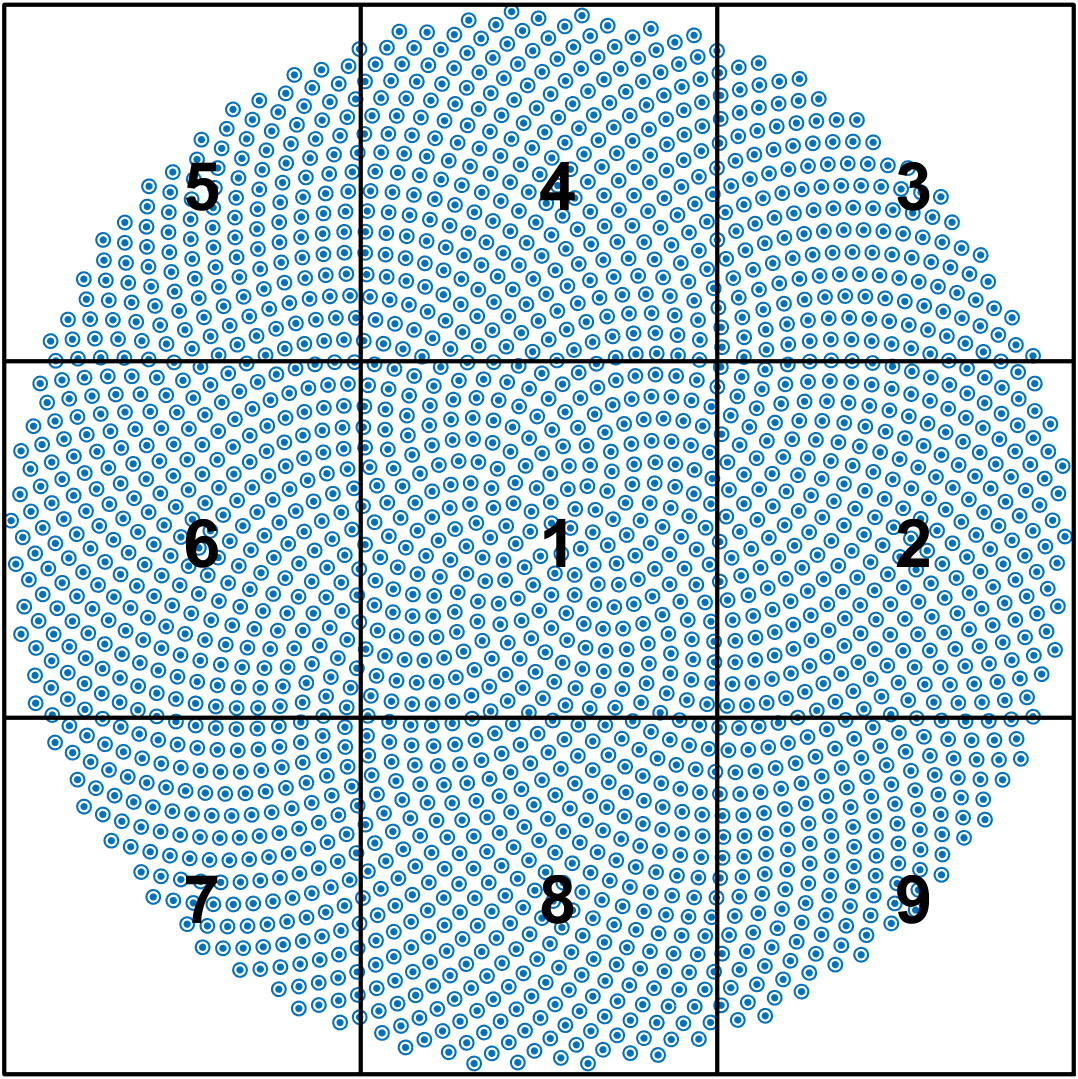}} \quad
    \subfloat[][]
    {\includegraphics[width=0.45\linewidth]{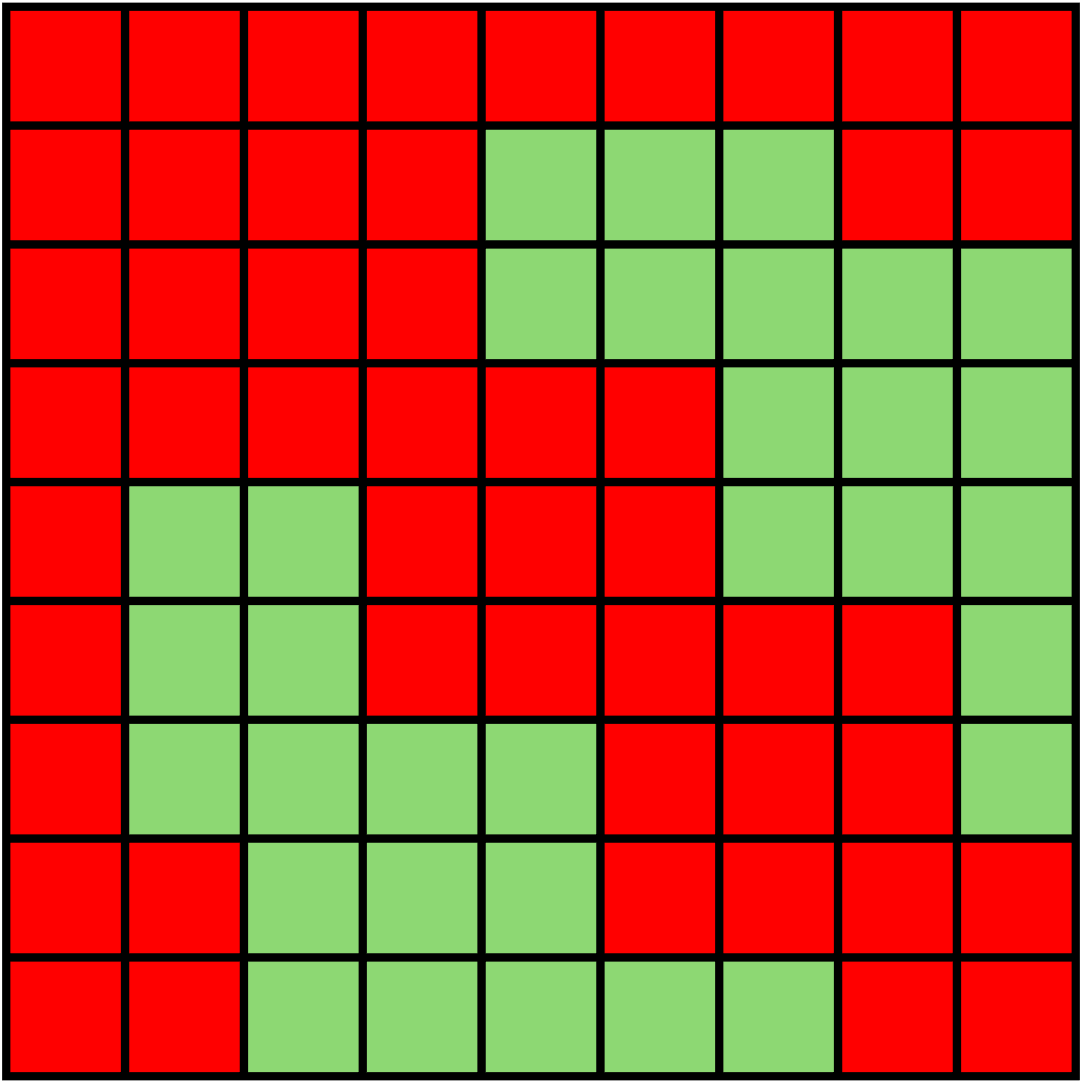}}
    \caption{Left: the meta-lens together with the grid (black) defining the blocks at the first grid level. Right: matrix $\ZZ$ partitioned according to the 9 blocks of meta-atom shown at left. Far distance interactions are marked in green, while near distance interactions are marked in red.}
    \label{fig_grid1}
\end{figure}

The classification in terms of near and far distance interactions is a relative concept related to a prescribed grid. Specifically, given a regular grid (see Figure \ref{fig_grid1}) superimposed to the metasurface, the interaction between meta-atoms that belong to the same block or a near block are classified as near, while interactions between non-adjacent blocks are classified as far.

\begin{definition}
Two blocks of the grid are near each other if they share at least one node.   
\end{definition}

For instance, matrix $\ZZ$ when partitioned block-wise according to the grid of Figure \ref{fig_grid1}, presents the near interaction blocks in red and the far interaction blocks in green. Matrix $\ZZ$ can be written as
\begin{equation}
\label{eq_split01}
    \ZZ = \ZZ_F^{(1)} + \ZZ_N^{(1)},
\end{equation}
where $\ZZ_F^{(1)}$ is the matrix made by only the far interactions (green block), while $\ZZ_N^{(1)}$ consists of only the near distance interactions (red blocks). The superscript highlights that the far-near splitting is at the first grid level. In formal terms, it results that
\begin{equation}
\label{eq_split02}
    \ZZ=\mathcal{F}_F^{(1)} \ZZ + \mathcal{F}_N^{(1)} \ZZ,
\end{equation}
where $\mathcal{F}_F^{(1)}$ is the linear operator selecting only the far distance interactions 
\[
\left( \mathcal{F}_F^{(1)}: \ZZ \mapsto \ZZ_F^{(1)} \right),
\]
while $\mathcal{F}_N^{(1)}$ selects only the near distance interactions 
\[
\left( \mathcal{F}_N^{(1)}: \ZZ \mapsto \ZZ_N^{(1)} \right).
\]

\subsection{The low-rank QR-decomposition}\label{sec_qrdec}
The far-near distance splitting is very convenient because each individual block accounting for far distance interactions (green blocks in Figure \ref{fig_grid1}) has a rank significantly lower than the size of the block itself, and it can be efficiently compressed via a low-rank QR-decomposition. Moreover, this rank decreases when the distance between the blocks increases.

Let $\left( \ZZ_F^{(1)} \right)_{i_1,j_1}$ be the $i_1j_1-$th block of matrix $\ZZ_F^{(1)}$. This block accounts for the interactions between the meta-atoms in block $i_1$ and the meta-atoms in block $j_1$, where $i_1$ and $j_1$ are far distance blocks. Matrix $\left( \ZZ_F^{(1)} \right)_{i_1,j_1}$ can be  accurately approximated within a prescribed tolerance $\varepsilon$, using a low-rank QR decomposition:
\begin{eqnarray}
\label{eq_02}
 \left( \ZZ_F^{(1)} \right)_{i_1,j_1}  \approx  \mathbf{Q}^{(1)}_{i_1,j_1} \mathbf{R}^{(1)}_{i_1,j_1}.
\end{eqnarray}
In equation \eqref{eq_02} $\left( \ZZ_F^{(1)} \right)_{i_1,j_1}$ is a $m \times n$ matrix, $\mathbf{Q}^{(1)}_{i_1,j_1}$ is a $ m \times r $ matrix, and $\mathbf{R}^{(1)}_{i_1,j_1}$ is a $r \times n$ matrix. The integer $r$ represents the rank of the block-to-block interaction matrix $\left( \ZZ_F^{(1)} \right)_{i_1,j_1}$ corresponding to the prescribed accuracy.

To find the impact of the low-rank QR-approximation of \eqref{eq_02}, it is convenient to evaluate the main computational cost and the memory storage. Specifically, the evaluation of the matrix-by-vector product $\left( \ZZ_F^{(1)} \right)_{i_1,j_1} \, \xx^{(1)}_{i_1j_1}$ with the ordinary formula calls for a computational cost proportional to $C_{\text{dir}} = \mathcal{O} (m \cdot n)$. Similarly, the memory required to store $\left( \ZZ_F^{(1)} \right)_{i_1,j_1}$ is proportional to $C_{\text{dir}}$. On the contrary, the computation of the matrix-by-vector product via the low-rank QR approximation of \eqref{eq_02} can be significantly reduced by performing the matrix multiplications in the proper order, i.e. by evaluating $\mathbf{Q}^{(1)}_{i_1,j_1} \left( \mathbf{R}^{(1)}_{i_1,j_1} \, \xx^{(1)}_{i_1j_1} \right)$. By doing this, the computational cost for computing the matrix-by-vector products is proportional to the smaller quantity $C_{\text{app}} = \mathcal{O} \left((m + n) r \right)$ and, similarly, the memory required to store both $\QQ^{(1)}_{i_1,j_1}$ and $\RR^{(1)}_{i_1,j_1}$ is proportional to the same quantity $C_{\text{app}}$.
When $r\ll m,n$, as it is the case for the interactions between blocks of meta-atoms at a far distance, it results that $ C_{\text{app}} \ll C_{\text{dir}}$, i.e. $G^{(1)}_{i_1j_1} \gg 1$, where $G^{(1)}_{i_1j_1} = C_{\text{dir}}/C_{\text{app}}$ is the compression gain for matrix $\left( \ZZ_F^{(1)} \right)_{i_1,j_1}$.

\subsection{The recursive scheme}
From a general perspective, the far/near splitting of \eqref{eq_split01} can be carried out on a coarse grid as well as on a finer grid. The the far distance interactions are compressed at the best on coarse grids made by large blocks. On the contrary, the impact of the near distance interactions is reduced on a finer grid made by small blocks. Large blocks are ideal for compressing the far distance interactions, while small blocks are ideal for reducing the number of near distance interactions.

The solution to these conflicting needs comes from a multilevel grid, where each block is subdivided into two parts along each orthogonal direction when moving from one subdivision level to the next. At each level, all the far distance interactions that have not been accounted for at the previous (coarser) level are compressed. The remaining near distance interactions are processed at the next refinement level, where they are partially compressed. For instance, when moving from the first to the second grid level, only $\ZZ^{(1)}_N$ is processed and split in far and near parts, i.e.
\begin{equation}
\begin{split}
\label{eq_rec01}
    \ZZ & = \mathcal{F}_F^{(1)} \ZZ + \mathcal{F}_N^{(1)} \ZZ \\
    & = \mathcal{F}_F^{(1)} \ZZ + \mathcal{F}_F^{(2)} \mathcal{F}_N^{(1)} \ZZ + \mathcal{F}_N^{(2)} \mathcal{F}_N^{(1)} \ZZ. 
\end{split}
\end{equation}
Similarly, at the third grid level only the near term $\mathcal{F}_N^{(2)} \mathcal{F}_N^{(1)} \ZZ$ is split in far and near parts. The far distance interactions blocks of $\mathcal{F}_F^{(1)} \ZZ$ are compressed at the first grid level, while the blocks of $\mathcal{F}_F^{(2)} \mathcal{F}_N^{(1)} \ZZ$ are compressed at the second grid level and so forth. The near distance interactions are progressively reduced in number, until the last grid level 
where they are retained without any approximation.
\begin{figure}[h!]
    \centering
    \subfloat[][]
    {\includegraphics[width=0.45\linewidth]{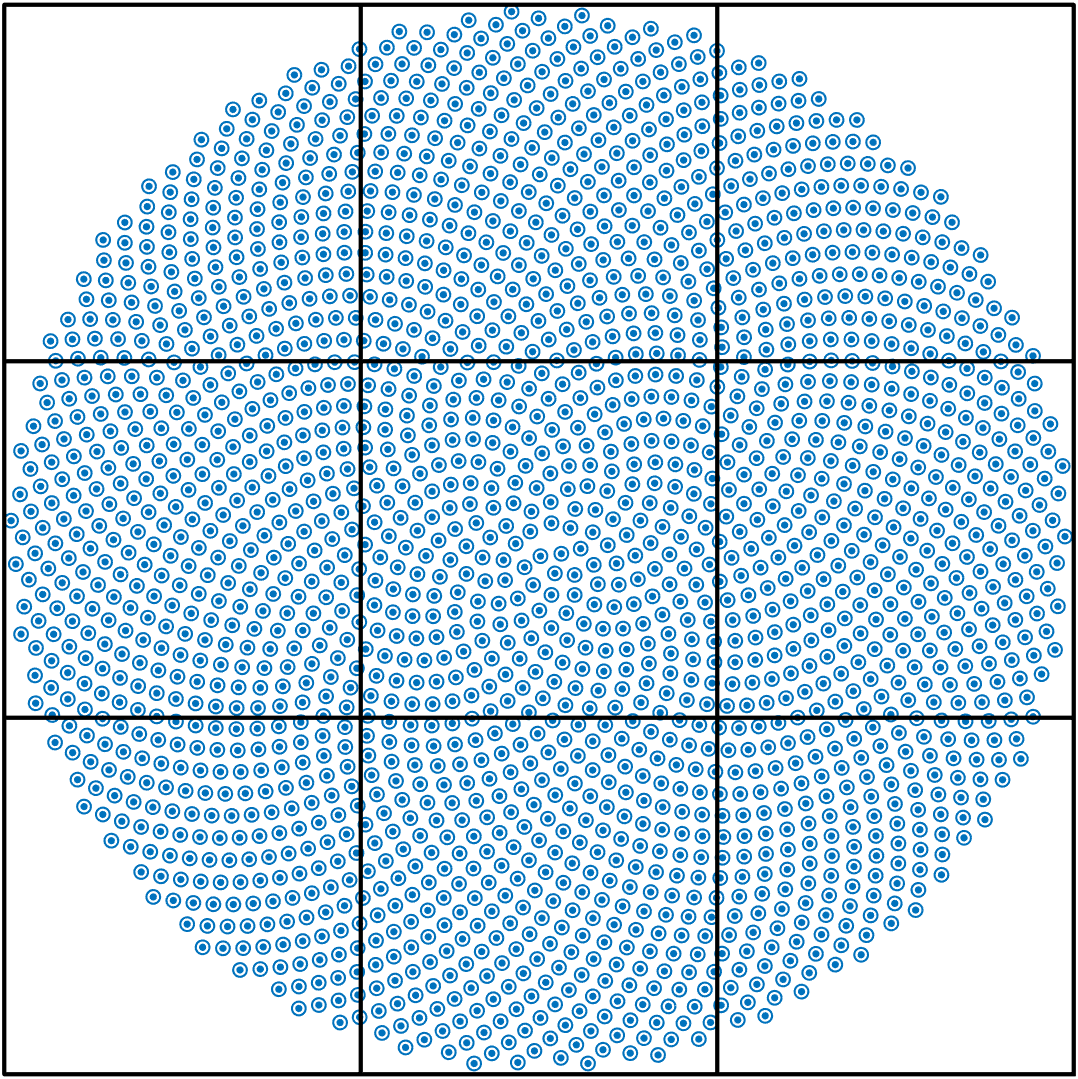}} \quad
    \subfloat[][]
    {\includegraphics[width=0.45\linewidth]{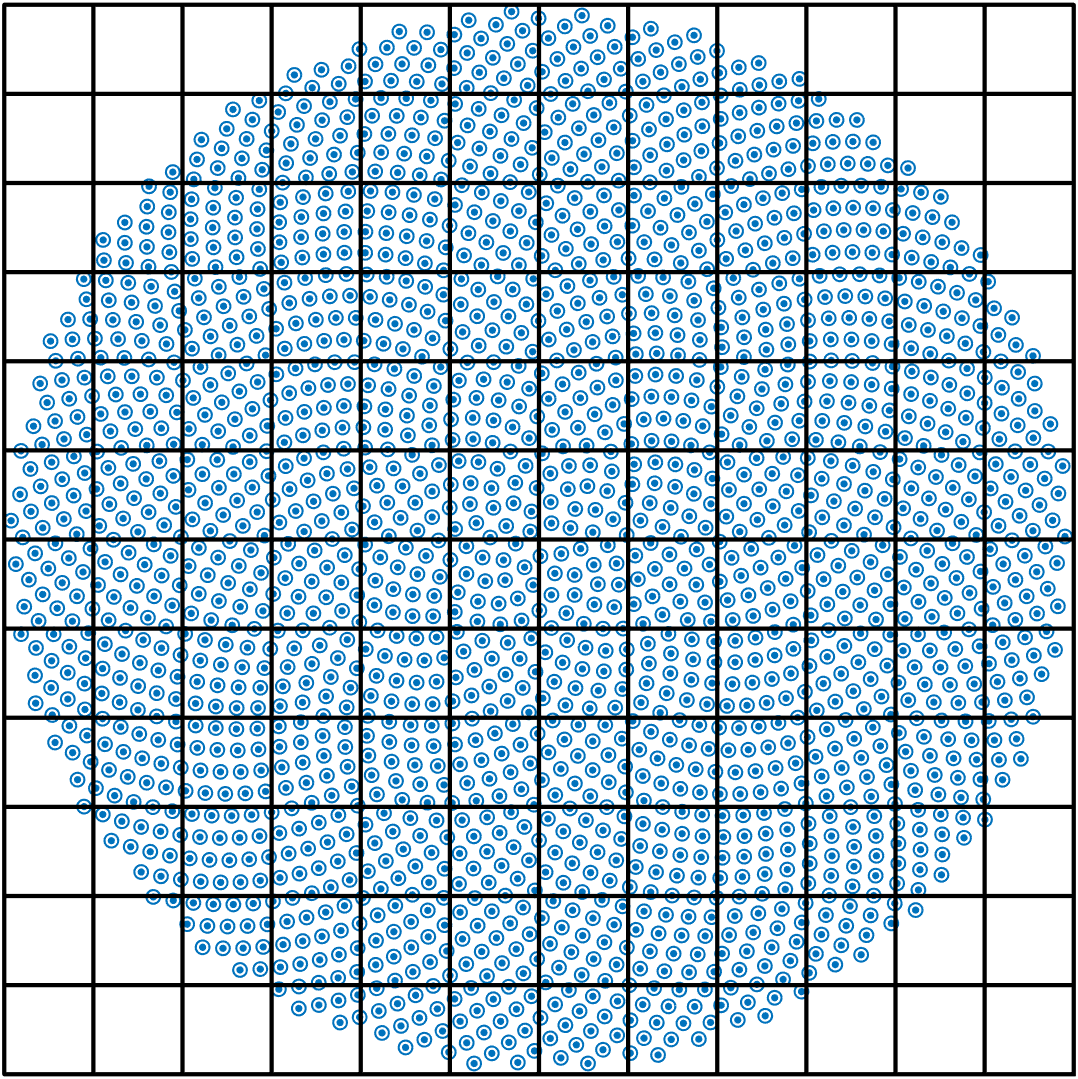}} 
    \caption{The grids for level 1 (left) and level 2 (right).}
    \label{fig_grid2}
\end{figure}
\begin{figure}[h!]
    \centering
    \includegraphics[width=0.6\linewidth]{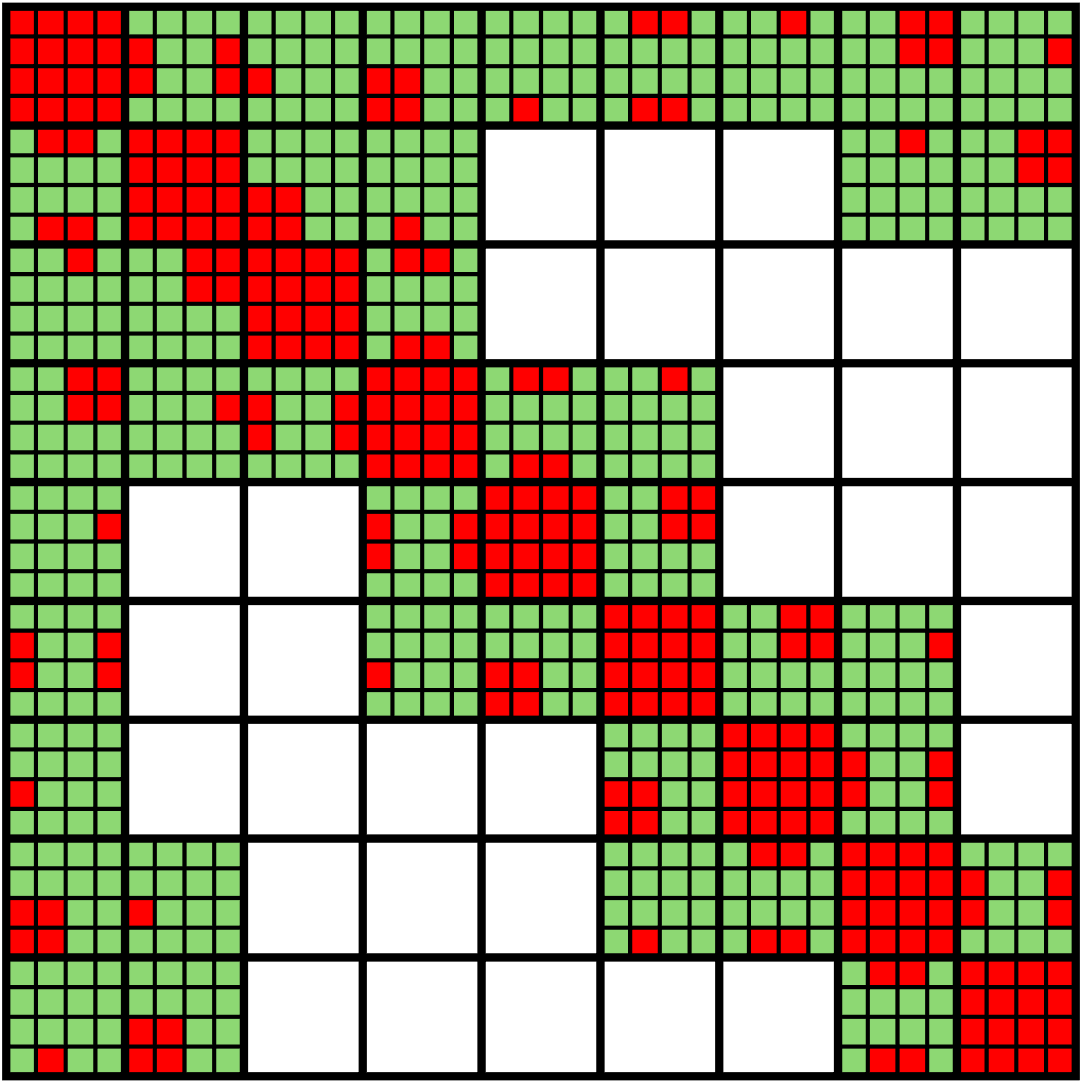}
    \caption{Matrix $\ZZ$ partitioned according to the 36  blocks at the second grid level. The new far distance interactions are marked in green, the current near distance interactions are marked in red, and the far distance interactions compressed at the previous level are marked in white.}
    \label{fig_grid3}
\end{figure}

The overall memory compression $G$ for the complete matrix $\ZZ$ is given by:
\begin{equation}\label{eqn_cgain}
G = \frac{N_{\text{dof}}^2}{N_{\text{QR}}},
\end{equation}
with:
\begin{equation}
    N_{\text{QR}} = N_{\text{far}} + N_{\text{near}} =
    \sum_{l=1}^{N_{\text{level}}}\sum_{i,j} \frac{m_i^{(l)}\times n_j^{(l)}}{G_{i,j}^{(l)}} + N N_p^2,
\end{equation}
where $N_{\text{QR}}$ is the size of the matrix after compression, $N_{\text{dof}}$ is the total degrees of freedom, $N_{level}$ is the number of grid levels, and $\left(\ZZ_F^{(l)}\right)_{i,j}$ has dimensions $m_i^{(l)}\times n_j^{(l)}$. The compression gain is inversely proportional to the required storage memory.

\subsection{The metasurface}
The proposed approach has been tailored for a metasurface made by an array of meta-atoms. Specifically, three ad-hoc customizations have been introduced and summarized below.

First, each meta-atom provides an individual unit where all of its DoFs are attributed to the same block of a given grid level. As discussed in \ref{app_PBG}, splitting a meta-atoms in smaller parts does not provide any benefit from the computational point of view. Therefore, each meta-atom, as a whole, is attributed to the block that contains its barycenter, at a given grid level.

Second, the finer grid level is chosen to have at most one meta-atoms per block to get $\ZZ_N$ with the smaller number of nonzero elements. In this manner, it results that $\ZZ_N=\ZZ_D$

Third, at the finer grid level, all the mutual interactions between distinct meta-atoms are compressed. This is because, thanks to the physical properties of the problem, the mutual interaction between distinct meta-atoms exhibits a rank deficiency.

\section{Numerical validation} \label{examples}
This section demonstrates the effectiveness of the proposed method through numerical validation, focusing on electromagnetic scattering problems involving electrically large metasurfaces. Following the framework discussed in Section \ref{nummod} and illustrated in Figure \ref{fig:qrbb}, we validate the proposed approach by analyzing computational time, solution accuracy, and memory usage as a function of the precision prescribed for the computation of the matrix-by-vector products.

The numerical validation is organized into several subsections, each addressing specific aspects of the method. After providing a detailed description of the test scenarios in Section \ref{sub_geo}, we first demonstrate the consistency of the proposed method. We show that both the memory requirements and solution accuracy scale with the precision prescribed for evaluating the matrix-by-vector products (Section \ref{sub_cons}). Subsequently, Section \ref{sub_perf} analyzes the overall performance of the approach, particularly highlighting the memory usage and the number of iterations required by the iterative solver, as the number of degrees of freedom increases. Finally, the efficiency of a parallel implementation is evaluated, with particular care paid to critical aspects such as load balancing across computational resources.

All numerical experiments were performed on a shared memory computing architecture (see details in Table 4). For the sake of consistency, each computation was performed using 64 cores (MPI processes).

The results presented throughout this section employ the GMRES method \cite{mybib:SS-GMRES} as iterative solver, paired with the preconditioner described in Section \ref{sec_prec}. The GMRES was selected because of its superior stability in the reduction of the residual, if compared to other iterative solvers, such as BICGSTAB. The maximum iteration count was set to \num{5000} without restarts, and convergence was numerically achieved for a residual below a threshold of $\num{1d-4}$.

\begin{table}[h!]
\centering
\captionsetup{justification=centering}
\caption{Machine used for computation}
\begin{adjustbox}{width=1\columnwidth,center}
\begin{tabular}{cc}
\hline\hline
\textbf{Specification}           & \textbf{Details}                       \\ \hline
Machine                          & HPE Proliant DL380 Gen10               \\ 
CPU                              & Intel(R) Xeon(R) Gold 6338 CPU @ 2.00GHz \\ 
Nodes                            & 2                                      \\ 
Cores per Node                   & 32                                     \\ 
L2 Cache                         & 80 MiB                                 \\ 
RAM                              & 2 TB                                   \\ \hline\hline
\end{tabular}\label{tab:machine_info}
\end{adjustbox}
\end{table}

\subsection{Description of the cases of interest}\label{sub_geo}

\begin{figure}[h!]
    \centering
    \includegraphics[width=0.95\linewidth]{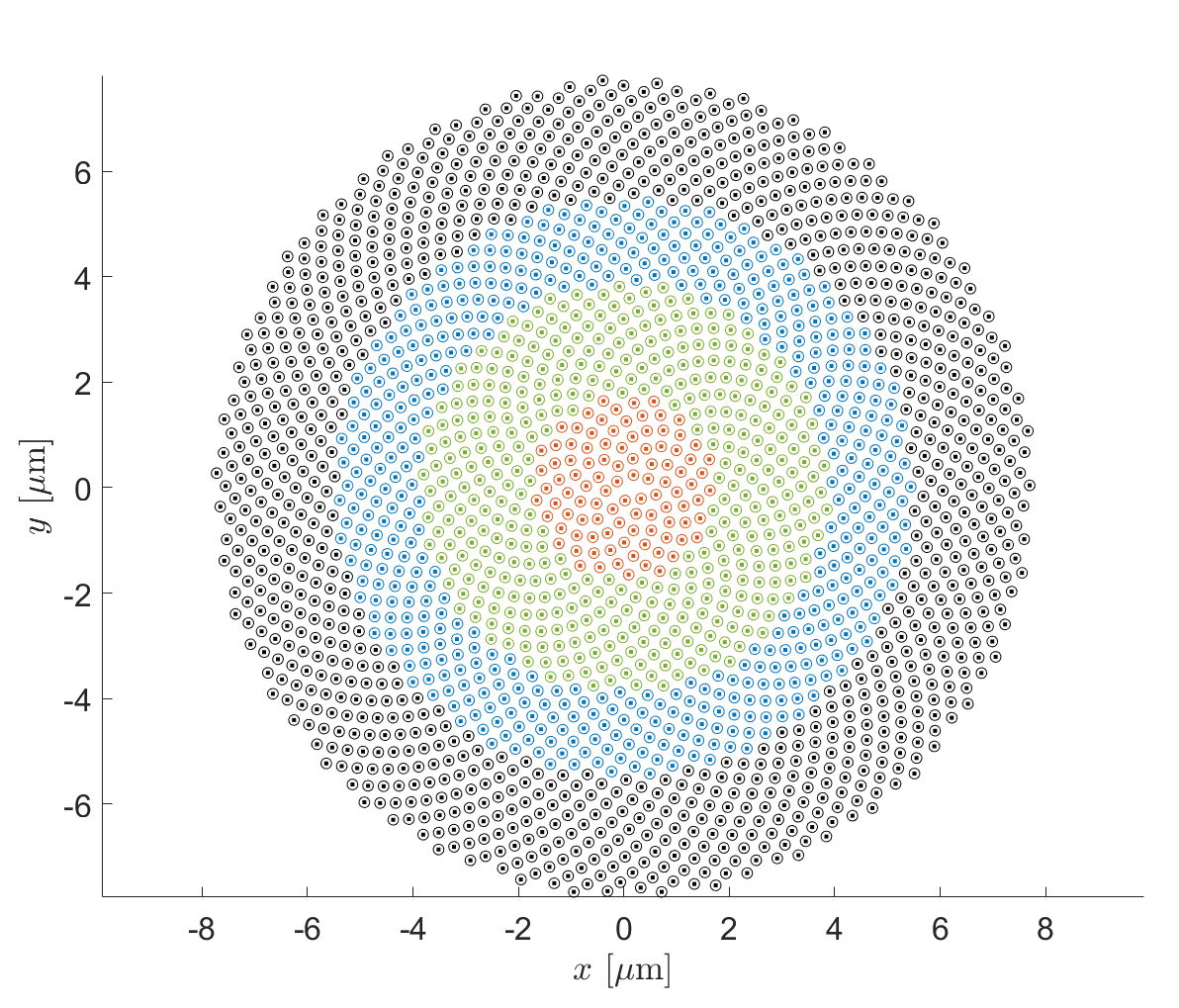}\\
    \caption{The four configurations considered in this paper. Each one is obtained by adding a certain number of meta-atoms to the previous one. Red $N=100$; green $N=500$; blue: $N=1000$; black: $N=2000$. The array radii are \qty{1.73}{\micro\meter}, \qty{3.87}{\micro\meter}, \qty{5.48}{\micro\meter}, \qty{7.75}{\micro\meter}, respectively}
    \label{fig:geo}
\end{figure}

The investigated metasurfaces, shown in Fig. \ref{fig:geo}, consist of arrays formed by replicating identical meta-atoms, of spherical shape. These spheres are arranged in an aperiodic structure known as Vogel's spiral \cite{Trevino2011}. This spiral arrangement is convenient for testing because individual meta-atoms can be added one by one without significantly altering the array density. Additionally, varying the number of meta-atoms modifies the array size, facilitating the investigation of the method's performance across different electrical lengths (ratios between the radius of the array and the free-space wavelength of the incident field).
The centers of the spheres for the Vogel's spiral are located at $(x_i, y_i, 0)=$ $r_i\left(\cos(i \alpha), \sin(i \alpha), 0 \right)$, for $i=1,\dots,N$, where $N$ is the number of meta-atoms, $r_i=c\sqrt{i}$, $c=R \sqrt{3}$ is a constant scaling factor, $R$ is the radius of each sphere, $\alpha=2\pi/\phi^2$ is the golden angle, and $\phi=(1+\sqrt{5}/2)$ is the golden ratio.
Four distinct arrays are analyzed, each composed of spheres with a radius of $R_s = \qty{100}{\nano\meter}$, arranged in configurations of $N = 100, 500, 1000$, and $2000$ meta-atoms. The smallest array ($N = 100$) is suitable for direct numerical analysis, serving as a reference case to validate the accuracy and consistency of the proposed numerical method by means of a direct solver, as detailed in the subsequent subsection. From the numerical perspective, all spheres share the same identical finite elements mesh. Table \ref{tab:case_studies} summarizes the main data for each array configuration.

The incident field is a plane wave at a frequency of $f = \qty{5d14}{\Hz}$, corresponding to a wavelength of \qty{600}{\nano\meter}. The spheres are made of gold \cite{johnson_optical_1972} characterized by a relative permittivity at the prescribed frequency $f$ of $\varepsilon_r = -9.428 + 1.513i$.

\begin{table}[h!]
\centering
\captionsetup{justification=centering}
\caption{Mesh data for the cases of interest}
\begin{adjustbox}{width=1\columnwidth,center}
\begin{tabular}{cccccc}
\hline \hline
\textbf{Meta-atoms} & \textbf{Points}& \textbf{Elements}  & \textbf{Loop} & \textbf{Star}  & \textbf{Total DOFs} \\ \hline
100  & 32100   & 25600  & 46500  & 9500 & 56000                                                          \\
500  & 160500  & 128000 & 232500 & 47500 & 280000           \\                                                     1000 & 321000  & 256000 & 465000 & 95000 & 560000                                                         \\
2000 & 642000  & 512000 & 930000 & 190000 & 1120000                                     \\          \hline \hline
\end{tabular}\label{tab:case_studies}
\end{adjustbox} 
\end{table}

\subsection{Consistency of the Method}\label{sub_cons}

As described in \Cref{nummod}, the effectiveness of the iterative method can be evaluated through three primary metrics: (i) the accuracy of the solution, (ii) the memory compression, and (iii) the reduction in the number of multiplications necessary to compute the matrix-by-vector product $\mathbf{Z I}$ \cite{10.1145/3126908.3126921}. Thus, it is crucial to analyze the trade-off between solution accuracy and memory compression, as a function of the prescribed accuracy in the computation of the products via compression. This latter quantity is termed $P_c$ and is defined as:
\begin{align*}
    P_c &= -\log(\varepsilon_{\text{prod}}),\\
    \varepsilon_{\text{prod}} &= \frac{\lVert \mathbf{Z}_{\text{QR}}\mathbf{u}-\mathbf{Z}\mathbf{u}\rVert}{\lVert \mathbf{Z}\mathbf{u}\rVert}, \\
    \mathbf{u} &= (1+j)[1,\dots,1]^T,
\end{align*}
where $\mathbf{Z}_{\text{QR}}$ represents the stiffness matrix after compression. The solution accuracy, denoted as $A_s$, is defined by means of the relative error with respect to the reference solution:
\begin{align*}
    A_s &= -\log(\varepsilon_{\text{sol}}),\\
    \varepsilon_{\text{sol}} &= \frac{\lVert \mathbf{J}_{\text{QR}} -\mathbf{J}_{\text{ref}} \rVert}{\lVert \mathbf{J}_{\text{ref}} \rVert},
\end{align*}
where $\mathbf{J}_{\text{ref}}$ represents the induced polarization density obtained via a direct method, i.e. LU decomposition, and $\mathbf{J}_{\text{QR}}$ denotes the solution achieved using the proposed iterative method method of Section \ref{sec_qr}.

\begin{figure}[h!]
    \centering
    \includegraphics[width=0.99\linewidth]{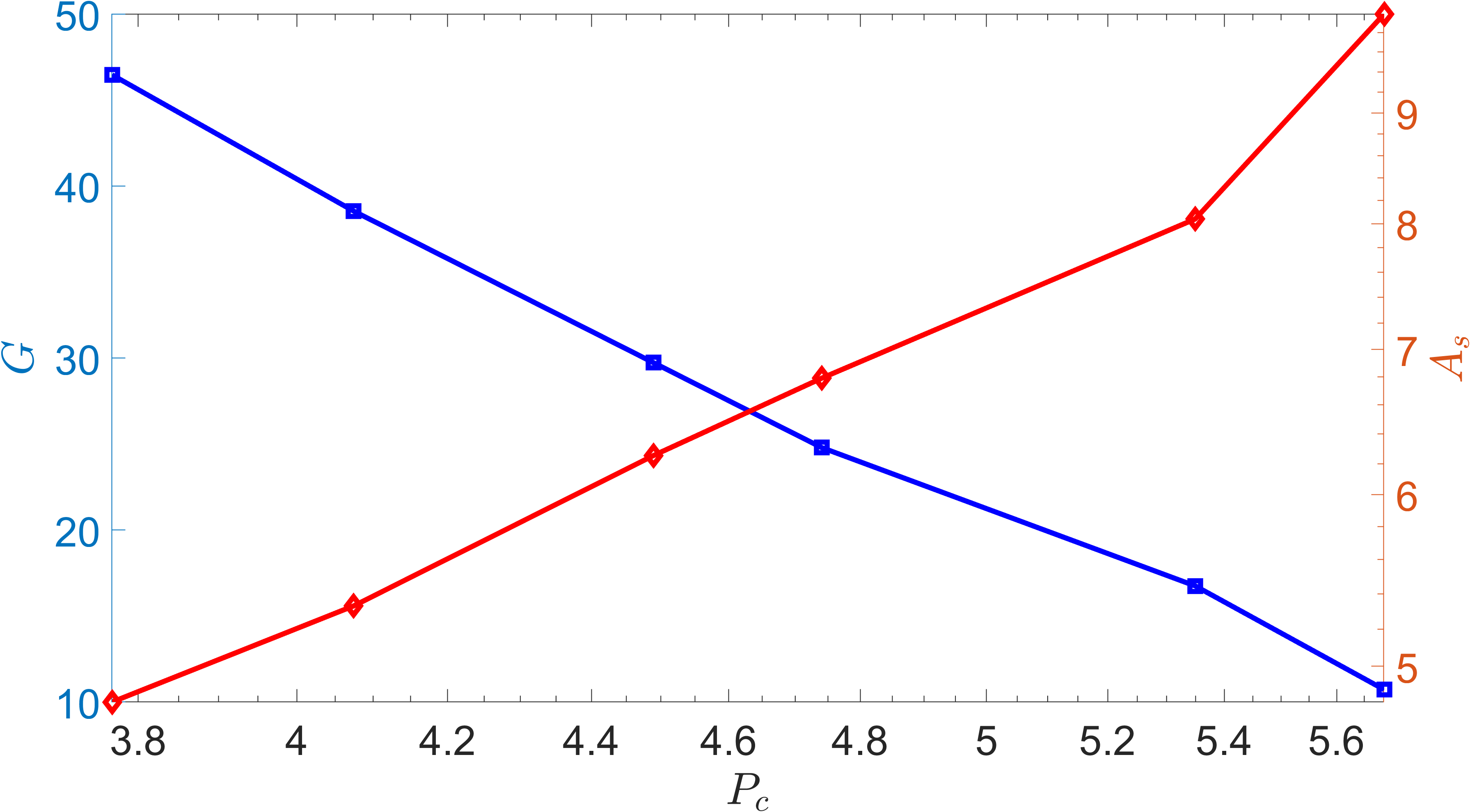}
    \caption{Method consistency:  Gain $G$ (blue line) and precision $A_s$ (red line) versus accuracy QR-recursive method.} 
    \label{fig:svd01}
\end{figure}

The numerical results for the reference case of 100 meta-atoms are shown in Figure \ref{fig:svd01}. As expected, increasing the input accuracy leads to a reduced compression gain and improved solution precision, clearly demonstrating the inverse relationship between these two metrics.
\subsection{Performances of the Method}\label{sub_perf}

The performances of the proposed method are evaluated by examining two primary metrics: the number of iterations required by the iterative solver (GMRES) to solve equation \eqref{eq_vie} and the compression gain as a function of the number of meta-atoms.

\begin{figure}[h!]
    \centering
    \includegraphics[width=0.99\linewidth]{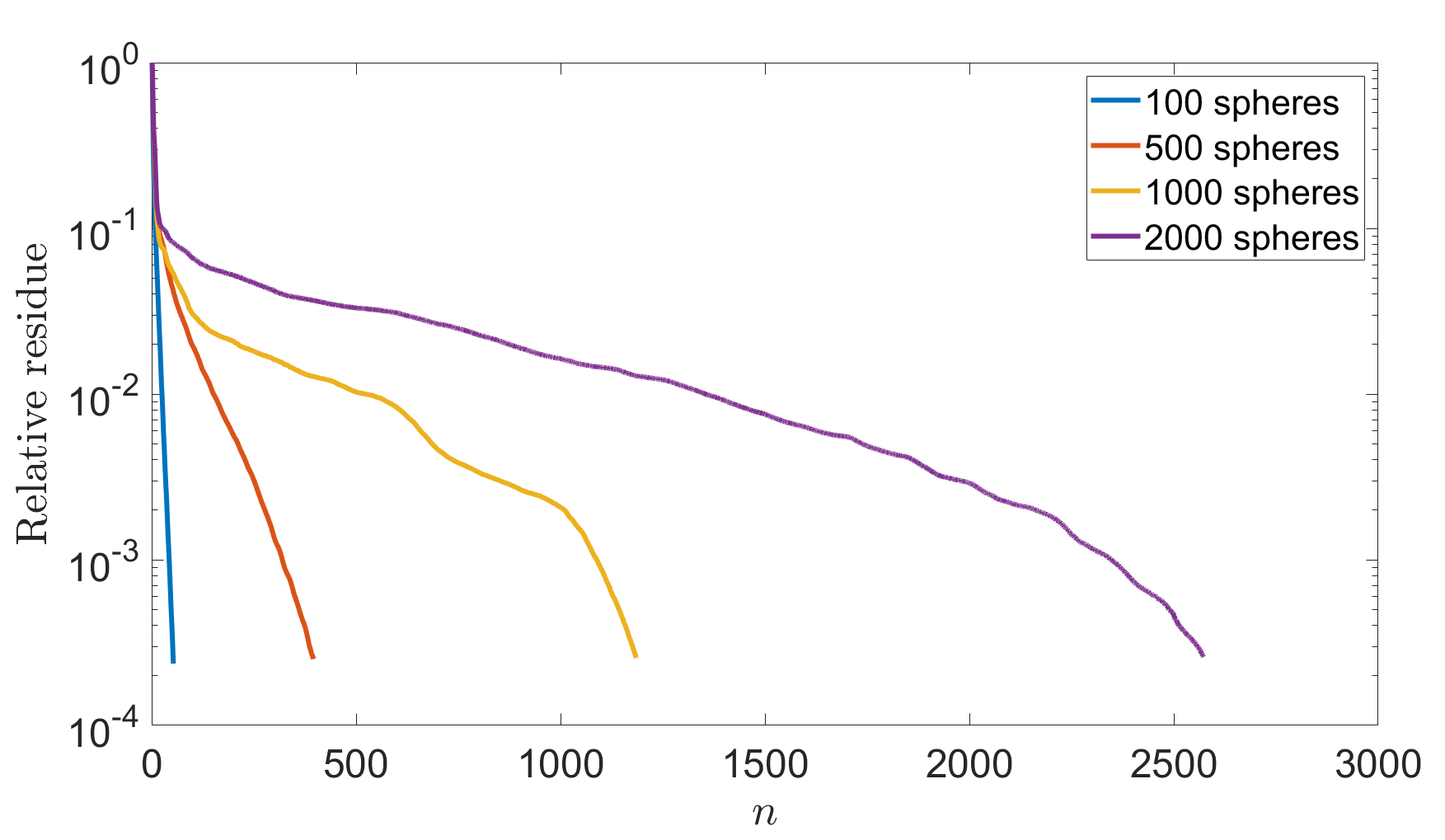}
    \caption{Relative residual norm versus the number of GMRES iterations for different meta-atom configurations.}
    \label{fig:it_gmres}
\end{figure}

Figure \ref{fig:it_gmres} presents the progression of the relative residual norm throughout the GMRES iterations, for arrays containing 100, 500, 1000, and 2000 meta-atoms.
As expected, the number of iterations grows when increasing the number of degrees of freedom (DoFs). Remarkably, even the largest scenario (2000 spheres) with approximately 1.1 million unknowns requires only 2600 iterations, clearly demonstrating the efficiency of the method. 
To further highlight the impact of the preconditioning strategy, we repeated the computation for the array containing 100 particles without using any preconditioner. In this case, the GMRES solver required 9067 iterations, compared to 53 iterations when the proposed preconditioner is employed.

Alternative preconditioners, such as incomplete LU decomposition or sparse approximate inverses \cite{benzi}, could potentially reduce iteration counts but at prohibitive costs in terms of both construction and storage, making them unsuitable for large and dense linear systems.

\begin{figure}[h!]
    \centering
    \includegraphics[width=0.99\linewidth]{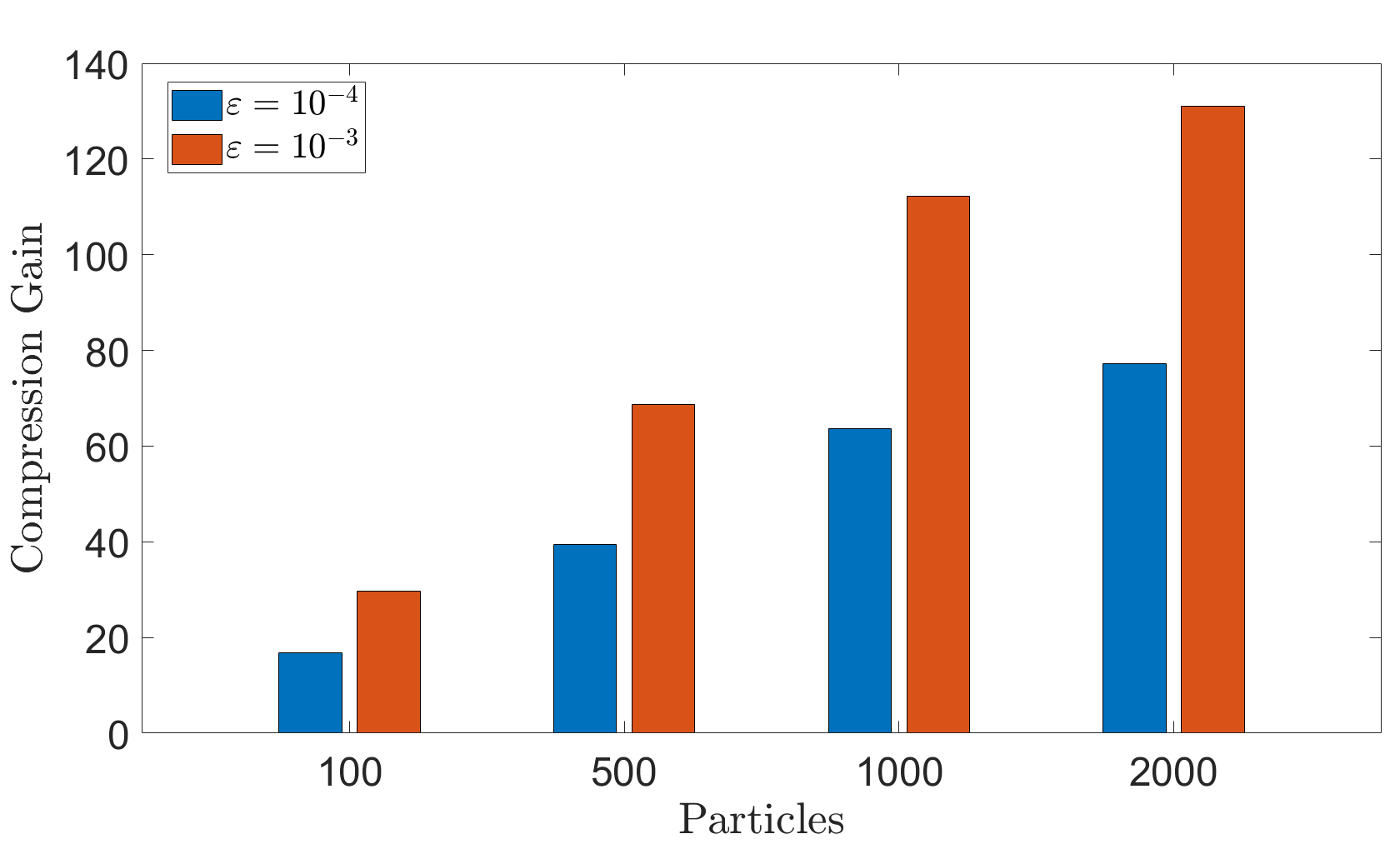}
    \caption{Compression gain as a function of the number of meta-atoms, at different values of QR tolerance $\varepsilon$ (see Section \ref{sec_qrdec}).}
    \label{fig:cg_par}
\end{figure}

Performance in memory usage is assessed through the compression gain $G$, defined in equation \eqref{eqn_cgain}, as a function of the number of unknowns. Figure \ref{fig:cg_par} illustrates how the compression gain improves with increasing meta-atom counts.

\begin{figure}[h!]
    \centering
    \includegraphics[width=0.99\linewidth]{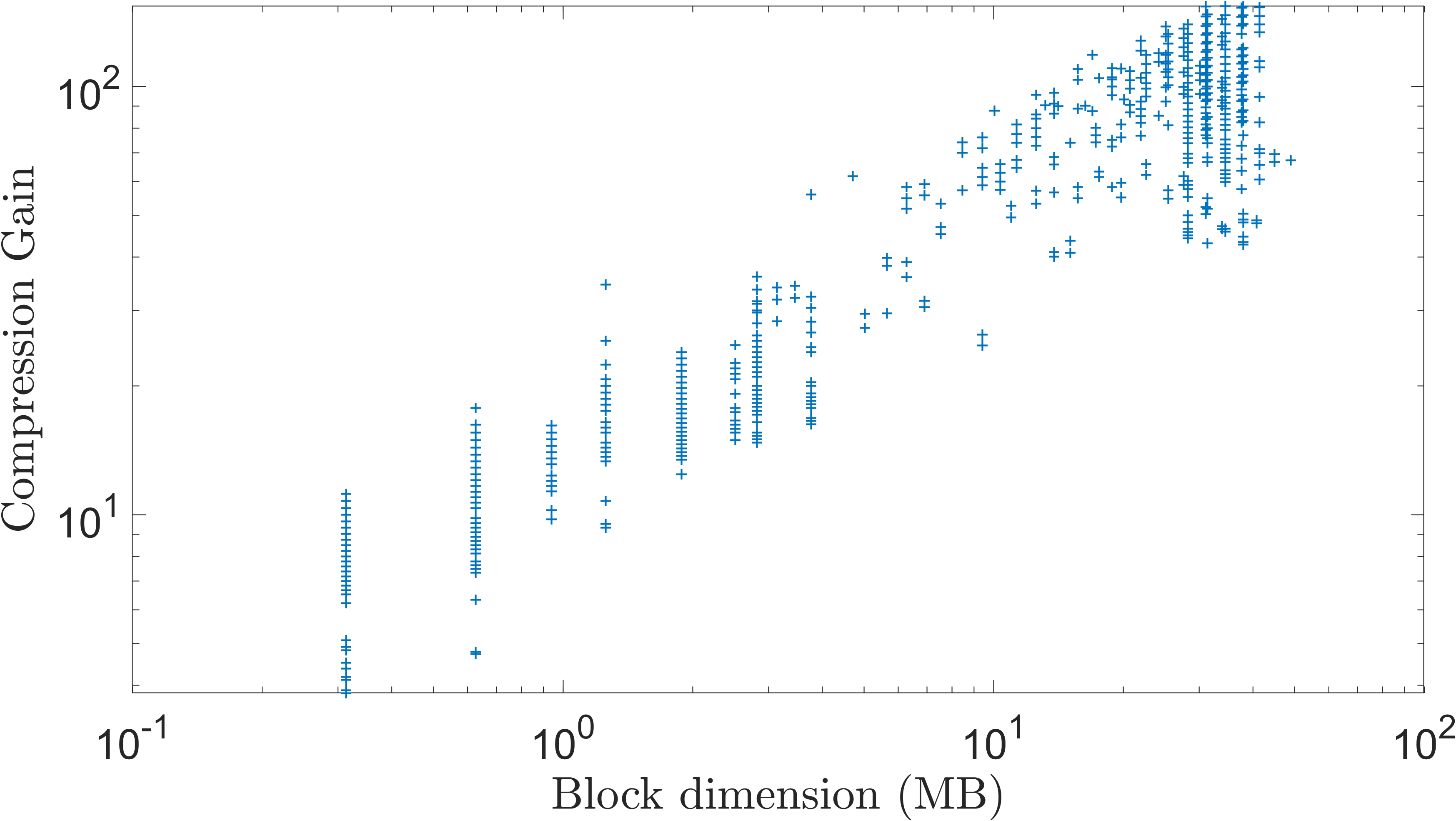}
    \caption{Compression gain as a function of the blocks size.}
    \label{fig:cg_block}
\end{figure}

To further characterize the compression behavior, Figure \ref{fig:cg_block} reports the compression gain of the individual interaction blocks arising in the recursive partitioning of the matrix $\mathbf{Z}$, shown as a function of the block size, i.e.  occupied memory. The results in Figure \ref{fig:cg_block} refer to the largest configuration considered in this work, namely the array composed of 2000 particles and QR tolerance $\varepsilon=10^{-3}$, for which the effect of the compression is most evident.

\subsection{Efficiency of the Parallelization Method}\label{sec_eff_p}
Three key factors determine the efficiency of the parallelization strategy: (i) balancing of the computational load during the assembly of the stiffness matrix, (ii) balancing of the memory required to evaluate the matrix-by-vector products, and (iii) balancing of the computation time required to evaluate the matrix-by-vector product. Computational load balancing is essential to ensure an even distribution of the computational workload, maximizing speedup in the construction of the compressed operator. Similarly, memory balancing plays a critical role, particularly in high-performance computing (HPC) environments where identical nodes are used. A well-balanced memory distribution prevents bottlenecks and allows efficient resource scaling as computational demand increases. Furthermore, balancing the computation time during the matrix-by-vector product is mandatory to achieve optimal speedup and enable efficient scalability as the number of processes increases. A suboptimal yet effective algorithm was implemented to distribute the workload efficiently, as detailed in \ref{appendix:Load}.  

Table \ref{tab:par_p} shows some statistics on the distribution of the computational load during the assembly of the matrix, in terms of the number of interactions assigned to each processor. Specifically, it includes the average number of interactions assigned to each processor \( \mu_a \), the standard deviation of the number of interactions per processor \( \sigma_a \), and the normalized standard deviation \( \sigma_a / \mu_a \). The low values of the normalized standard deviation indicate an excellent balance achieved during the assembly phase, ensuring an even distribution of computational efforts across processors.

\begin{table}[!h]
\centering
\captionsetup{justification=centering}
\caption{Mean value \( \mu_a \) of the number of interactions assigned to each processor during assembly, corresponding standard deviation \( \sigma_a \), and normalized standard deviation \( \sigma_a / \mu_a \).}
\begin{adjustbox}{width=1\columnwidth,center}
\begin{tabular}{ccccc}
\hline \hline
\textbf{Meta-atoms} & \textbf{Total interactions} & \( \mu_a \) & \( \sigma_a \) & \( \sigma_a / \mu_a \) \\ \hline
100  & 9900     & 155   & 0.47 & \num{3.00d-3}   \\
500  & 249500   & 3898  & 0.50 & \num{1.28d-4}      \\                                                     
1000 & 999000   & 15609 & 0.49 & \num{3.13d-5}   \\
2000 & 3986943  & 62296 & 0.13 & \num{2.00d-6} \\          \hline \hline
\end{tabular}\label{tab:par_p}
\end{adjustbox} 
\end{table}

Table \ref{tab:par_bal} shows similar statistics but related to the memory allocation per processor, including the average value \( \mu_m \), the standard deviation \( \sigma_m \), and the normalized standard deviation \( \sigma_m / \mu_m \). Even for this aspect, the consistently low values of the normalized standard deviation confirm that the memory distribution remains well-balanced, preventing potential bottlenecks in large-scale computations.

\begin{table}[!h]
\centering
\captionsetup{justification=centering}
\caption{Average allocated memory \( \mu_m \) per processor, corresponding standard deviation \( \sigma_m \), and normalized standard deviation \( \sigma_m / \mu_m \) for different test cases with varying numbers of meta-atoms.}
\begin{adjustbox}{width=1\columnwidth,center}
\begin{tabular}{ccccc}
\hline \hline
\textbf{Meta-atoms} & \textbf{Total memory (GB)} & \( \mu_m \) (GB) & \( \sigma_m \) (GB) & \( \sigma_m / \mu_m \) \\ \hline
100  & 1.57     & 0.024   & \num{4.90d-3} & 0.198   \\
500  & 16.99   &  0.265 & \num{1.05d-2} & 0.040      \\                                                     
1000 & 41.66   & 0.651 & \num{1.84d-2} & 0.028   \\
2000 & 142.67  & 2.229 & \num{2.92d-2} & 0.013 \\          \hline \hline
\end{tabular}\label{tab:par_bal}
\end{adjustbox} 
\end{table}

\subsection{Timing Information}

\begin{table}[!h]
\centering
\captionsetup{justification=centering}
\caption{Assembly Time for various array sizes.}
\begin{adjustbox}{width=0.7\columnwidth,center}
\begin{tabular}{ccc}
\hline \hline
\textbf{Array Size} & \textbf{Assembly (s)} & \textbf{Compression (s)} \\ \hline 
100           & 727                       & 70                         \\ 
500           & 11579                     & 164                        \\ 
1000          & 41012                     & 11568                      \\ 
2000          & 116211                    & 65369                      \\ \hline\hline
\end{tabular}\label{tab:performance0}
\end{adjustbox} 
\end{table}

The construction process of the compressed matrix-by-vector operator involves two main steps: (i) local matrix assembly, and (ii) subsequent compression employing the full-cross approximation technique \cite{Kurz2002}. 
First, the assembly times required to construct the operator are shown in Table \ref{tab:performance0} for the scenarios considered. It is evident that, as the size of the array increases, the compression time becomes a significant fraction of the overall duration of the assembly. 
\begin{table*}[h!]
    \centering
    \captionsetup{justification=centering}
    \caption{Performance analysis of different stages of the GMRES algorithm across multiple problem sizes (100 to 2000). Values marked with an asterisk ($^*$) indicate extrapolated estimates. All reported timings are expressed in seconds.}
    \begin{adjustbox}{width=0.88\textwidth,center}
    \begin{tabular}{clccccccc}
    \hline \hline
        \textbf{Case} & & \textbf{100} & \textbf{500} & \textbf{1000} & \textbf{2000} \\
        \hline 
        $N_\text{prd}$& number of matrix products ($\mathbf{Z}\,\mathbf{x}$) & 107 & 793 & 2371 & 5151 \\
        $N_\text{Pcnd}$ & number of  preconditioner applications  & 54 & 397 & 1186 & 2576 \\
         ${\rho}_\text{prd}$ & $\mathbf{Z} \mathbf{x}$ product speed-up & \textbf{13.98}  & \textbf{53.52}  & $\textbf{95}^*$ & $\textbf{96}^*$ \\
        $t_{prd}^{\text{QR}}$ & overall QR product $\mathbf{Z}_{QR} \mathbf{x}$ time  & 1.15 & 77.6 & 530 & \num{4.47d3} \\
        $t_{prd}^{\text{ord}}$ & overall ordinary $\mathbf{Z} \mathbf{x}$  product-time & 16.03 & \num{4.15d3} & {\num{4.97d4}}$\,^*$ & {\num{4.31d5}}$\,^*$ \\
        $t_\text{Pcnd}$ &  overall preconditioner application time & \num{2.42d-1} & 11.68 & 65 & 282 \\
        $t_\text{Kryl}$ & Krylov algebraic operations Time & 17.09 & 259 & \num{4.12d3} & \num{5.68d4} \\
        $t^{\text{QR}}_{tot}$ & QR-GMRES overall exc. time & 18.48 & 349 & \num{4.71d3} & \num{6.16d4} \\
        $t^\text{ord}_{tot}$ & Ordinary GMRES overall exc. time & 33.36 & \num{4.42d3} & {\num{5.38d4}$^*$} & \textbf{\num{4.88d5}}$^*$ \\
        ${\rho}_\text{GMRES}$ & \textbf{overall GMRES speed-up}  & \textbf{1.8} & \textbf{12.68} & $\textbf{11.42}^*$ & $\textbf{7.93}^*$ \\
        \hline \hline
    \end{tabular}
    \end{adjustbox}
    \label{tab:transposed_table}
\end{table*}

Second, Table \ref{tab:transposed_table} provides an in-depth analysis of inversion timing, focusing on the most computationally demanding operations: (i) matrix-by-vector product ($\mathbf{Z}\mathbf{x}$), (ii) preconditioner application ($\mathbf{P}^{-1}\mathbf{y}$), and (iii) Krylov subspace algebraic manipulations necessary for obtaining the minimal solution. We implemented a loop employing the Arnoldi process to minimize the residual error iteratively. It is worth noting that an inefficient preconditioner results in a quadratic growth of the dimension of the Krylov subspace, thus significantly slowing the convergence. To assess the efficiency of the proposed approach, Table \ref{tab:transposed_table} compares the GMRES inversion with and without the QR-recursive compression. The reference implementation without QR-recursive compression, referred to \textit{ordinary} implementation, utilizes optimized PBLAS routines to ensure peak efficiency. The presented results clearly illustrate substantial performance gains due to the proposed QR-recursive compression technique. Specifically, the ratio $\rho_{prd}$ of the computation times between ordinary and compressed matrix-by-vector products during the GMRES inversion, demonstrates approximately a hundred-fold acceleration. However, this significant compression translates to an overall GMRES execution speed-up of approximately tenfold compared to the standard GMRES implementation. Thus, the total performance gain does not scale linearly with the compression improvement observed in the matrix products. This discrepancy arises primarily due to limitations in the efficiency of the preconditioner. Despite this limitation, expanding the Krylov subspace remains indispensable for ensuring high solution accuracy.

\section{Conclusion} \label{conclusions}
We proposed a fast and robust integral equation solver specifically designed for efficiently handling electrically large particle arrays, extending to linear sizes of many (dozens) wavelengths, a scenario frequently encountered in electromagnetic modeling of metasurfaces and metalenses. Our method integrates a novel QR-recursive compression technique within a volume integral equation (VIE) framework, leveraging an iterative solver coupled with a preconditioner explicitly tailored to exploit the geometric structure of the array. This combination significantly accelerates the iterative solver, maintaining high accuracy.
This represents a first and key contribution in view of the simulation of suitable large-scale electromagnetic problems where the scatterer extends for thousands of wavelengths.

Our numerical experiments validated the consistency and scalability of the proposed approach.

We identified and analyzed the three critical factors that dictate the efficiency of our parallelization strategy: (i) balancing of the computational load during the assembly of the stiffness matrix, (ii) balancing of the memory required to evaluate the matrix-by-vector products, and (iii) balancing of the computational time required during the matrix-by-vector products.

By analyzing the scattering from finite arrays of spheres and performing comparative studies against conventional iterative VIE solvers without compression, we quantitatively demonstrated that the proposed QR-recursive compression, combined with a proper pre-conditioner and the GMRES iterative method, achieves approximately a tenfold reduction in computational time, plus substantial memory savings. Our extensive validation further confirmed the solver's effectiveness even for scenarios involving millions of degrees of freedom, where previously employed VIE approaches are impractical.

These findings highlight that our proposed solver represents a highly valuable computational tool, uniquely positioned to significantly advance the analysis and design capabilities in the field of metasurfaces and complex electromagnetic systems.

\section{Acknowledgment}
This work was supported by the Italian Ministry of University and Research under the PRIN-2022, Grant Number 2022Y53F3X \lq\lq Inverse Design of High-Performance Large-Scale Metalenses\rq\rq.

\appendix

\section{Properties of $\J$}
This Appendix reviews the main properties of the unknown vector field $\J$ as defined in \eqref{eq:const_rel} and introduced in \cite{Rubinacci20062977,Miano20102920,Forestiere20171224}.

\begin{figure}[h!]
    \centering
    \subfloat[][]
    {\includegraphics[width=0.65\columnwidth]{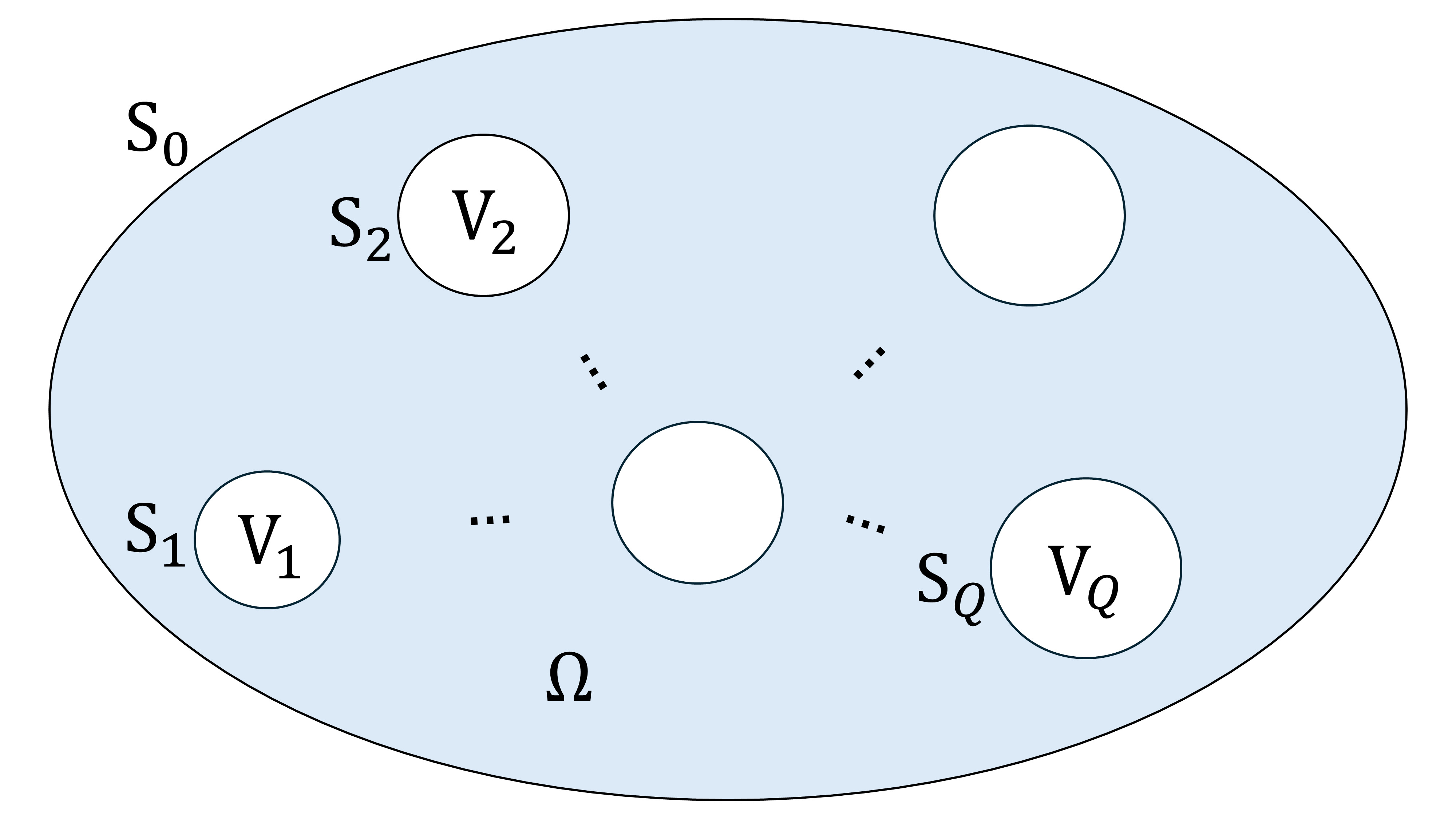}} \\
    \subfloat[][]
    {\includegraphics[width=0.65\columnwidth]{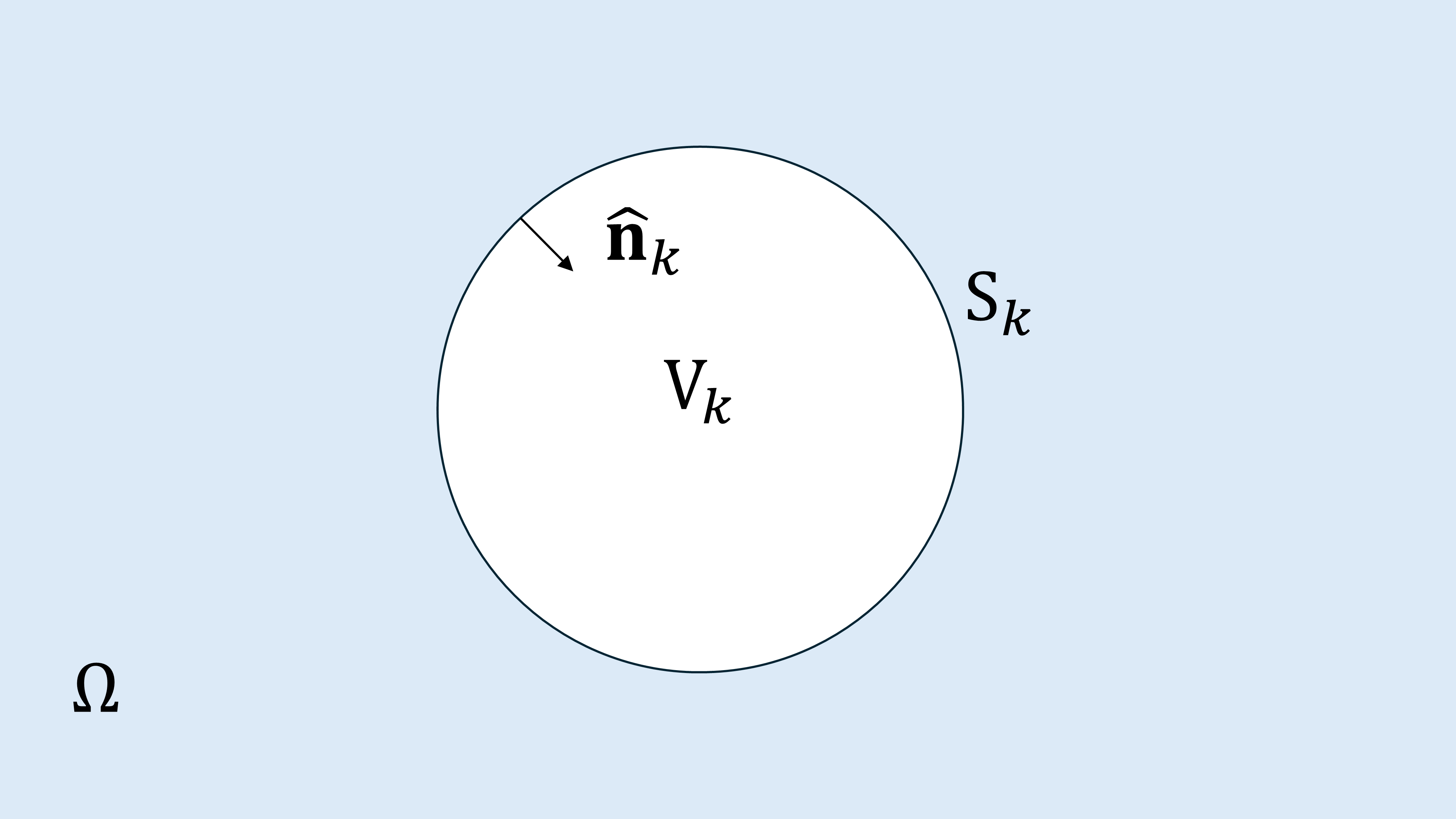}} 
    \caption{Domain $\Omega$ and its internal cavities $V_1, V_2, \dots, V_Q$ (top) and particular of the cavity $V_k$ with its boundary $S_k$ and the normal unit vector $\hat{\mathbf{n}}_k$ (bottom).}
    \label{Geo_Meta}
\end{figure}

Without loss of generality, let's $\Omega$ be made by one meta-atom, only. Aside from being a connected domain, this meta-atom may have an arbitrary topology. For example, the boundary $\partial \Omega$ is assumed to be the union of one or more distinct connected components $S_0,\ldots, S_Q$, where $S_0$ is the outer boundary while $S_1,\ldots,S_Q$ are the boundaries of cavities internal to $\Omega$ (see Figure \ref{Geo_Meta}(a)). The volumetric region occupied by each cavity and bounded by $S_k$ is denoted as $V_k$.

The main assumptions are that (i) the material properties $\sigma$ and $\chi$ are constant in each meta-atom, (ii) that each meta-atom is free from source charges in $\Omega$, and (iii) its cavities $V_1,\ldots,V_Q$, if any, are free from source charges. These assumptions are met in the majority of meta-surfaces occurring in practice.

First, $\J$ is divergence-free. Indeed, it results that
\begin{align}
\nabla \cdot \J & = (\sigma + j \omega \chi ) \nabla \cdot \E \\
 & = \frac{\sigma + j \omega \chi}{\varepsilon_0 (1+\chi) } \nabla \cdot \DD \\
 & = 0,
\end{align}
where the first two equalities come from the material properties $\sigma$ and $\chi$ being constants, while the latter equality comes from the absence of source charge in $\Omega$.

On the boundary of cavity $V_k$ (see Figure \ref{Geo_Meta}(b)), the unpaired surface charge density $ \rho^U_S$ is responsible for a jump of the normal component of the conduction current $\sigma \E$ and of the electric displacement $\DD$ (see \cite{HausMelcher1989}):
\begin{align}
    \label{eq001}    \sigma \E_- \cdot \mathbf{\hat{n}}_k & = j \omega \rho^U_S \text{ on } S_k\\
    \label{eq000}    \DD_+ \mathbf{\hat{n}}_k - \DD_- \mathbf{\hat{n}}_k& =  \rho^U_S \text{ on } S_k,
\end{align}
where $\DD_\pm (P)= \lim_{h \to 0^+} \DD(P \pm h \mathbf{\hat{n}}_k)$ and $\E_\pm (P)= \lim_{h \to 0^+} \E(P \pm h \mathbf{\hat{n}}_k)$, for any point $P$ on $S_k$.

Second, neutrality of the surface charge is achieved on the boundary of each cavity:
\begin{equation}
    \int_{S_k} \J_- \cdot \mathbf{\hat{n}}_k \text{d} S = 0, \, \forall k=1,\ldots,Q.
\end{equation}
Indeed, in cavity $V_k$, from the divergence theorem, it results that:
\begin{align}
0 & = - \int_{V_k} \nabla \cdot \DD \, \text{d} \tau \\
    & = \int_{S_k}\DD_+ \cdot \mathbf{\hat{n}}_k \text{d} S \\
    & = \int_{S_k} \rho_S^U \text{d} S + \int_{S_k}\DD_- \cdot \mathbf{\hat{n}}_k \text{d} S \\
    & = \frac{\sigma + j \omega \varepsilon}{j \omega} \int_{S_k} \E_- \cdot \mathbf{\hat{n}}_k \text{d} S \\
    & = \frac{\sigma + j \omega \varepsilon}{j \omega(\sigma + j \omega \varepsilon_0 \chi)} \int_{S_k} \J_- \cdot \mathbf{\hat{n}}_k \text{d} S,
\end{align}
where the first equality comes from the assumption that $V_k$ is free from source charge, the second equality from the divergence theorem, the third equality from \eqref{eq000}, the fourth equality from \eqref{eq001} and the constitutive relationship $\DD = \varepsilon \E$, and the last one from \eqref{eq:const_rel}.

Third, charge neutrality is achieved on $S_0$, the outer boundary of $\Omega$. Indeed, the divergence theorem yields
\begin{align}
    0 & = \int_{\Omega} \nabla \cdot \J \text{d} \tau \\
    & = \int_{S_k} \J_- \cdot \mathbf{\hat{n}}_k \text{d} S + \sum_{k=1}^Q \int_{S_k} \J_- \cdot \mathbf{\hat{n}}_k \text{d} S \\
    & = \int_{S_k} \J_- \cdot \mathbf{\hat{n}}_k \text{d} S.
\end{align}

\section{(Non) Splitting the meta-atoms}
\label{app_PBG}

This Appendix shows that the splitting into different blocks of the DoFs of an individual meta-atom is not convenient from the computational point of view.

\begin{figure}[h!]
    \centering
    \includegraphics[width=0.99\linewidth]{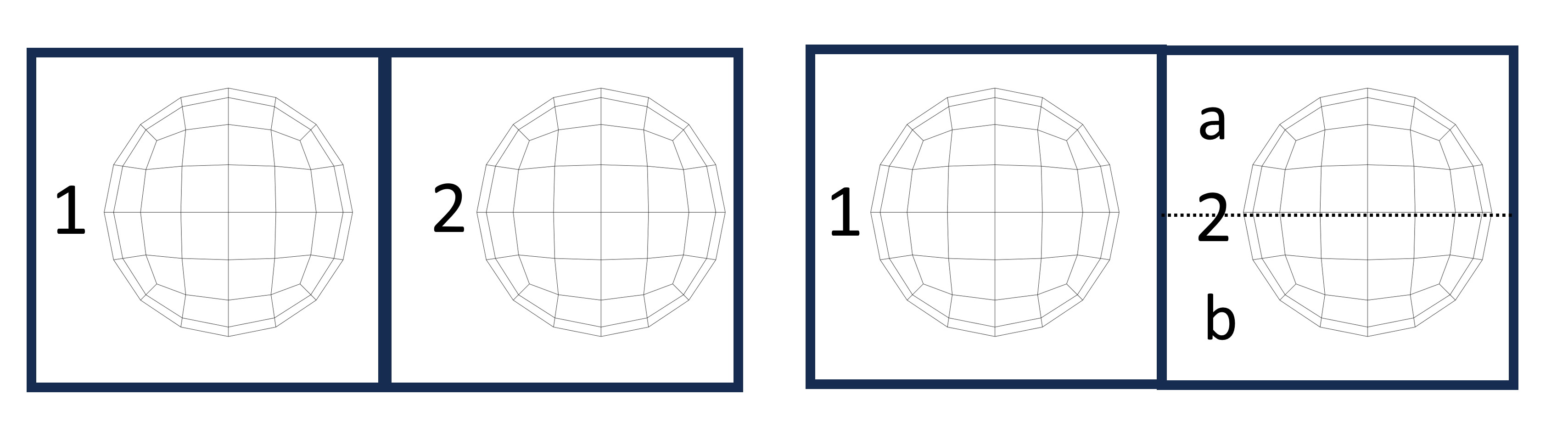}
    \caption{Two approaches for computing interaction matrix between two meta-atoms. Left: each meta-atom DOFs are assigned to a single block (1 and 2) : resulting in one matrix interaction. Right: DoFs related to the right meta-atom (the  block 2)  are split between two subblock blocks (namely a and b): resulting in the sum of two matrices  interactions..}
\label{fig:compgainvsrelres0}
\end{figure}

For the sake of simplicity, it is considered a simple configuration involving only two meta-atoms. Specifically, the performances regarding the computation of the mutual interactions between these two meta-atoms are evaluated in two scenarios. In the first one (Figure \ref{fig:compgainvsrelres0}(left)), the individual DoFs of two meta-atoms are not split into different blocks. In the second one (Figure \ref{fig:compgainvsrelres0}(right)), the individual DoFs of one meta-atom are split into two different blocks.

In the first case, the mutual interaction between the meta-atoms is approximated and compressed via the QR decomposition as
\begin{eqnarray}
\label{eq_singlebox}
 \ZZ_{1,2}  \approx  \mathbf{Q}_{1,2}\mathbf{R}_{1,2},
\end{eqnarray}
where the subscripts refer to the meta-atoms. In the second case, it turns out that
    \begin{align}
\label{eq_splitbox}
 \ZZ_{1,2} & = \ZZ_{1,2a} + \ZZ_{1,2b}   \\
 & \approx \mathbf{Q}_{1,2a}\mathbf{R}_{1,2a} + \mathbf{Q}_{1,2b}\mathbf{R}_{1,2b},
\end{align} 
where $2a$ $(2b)$ refers to the DoFs of meta-atom number 2 falling in block $a$ $(b)$.

\begin{figure}[h!]
    \centering
    \includegraphics[width=0.99\linewidth]{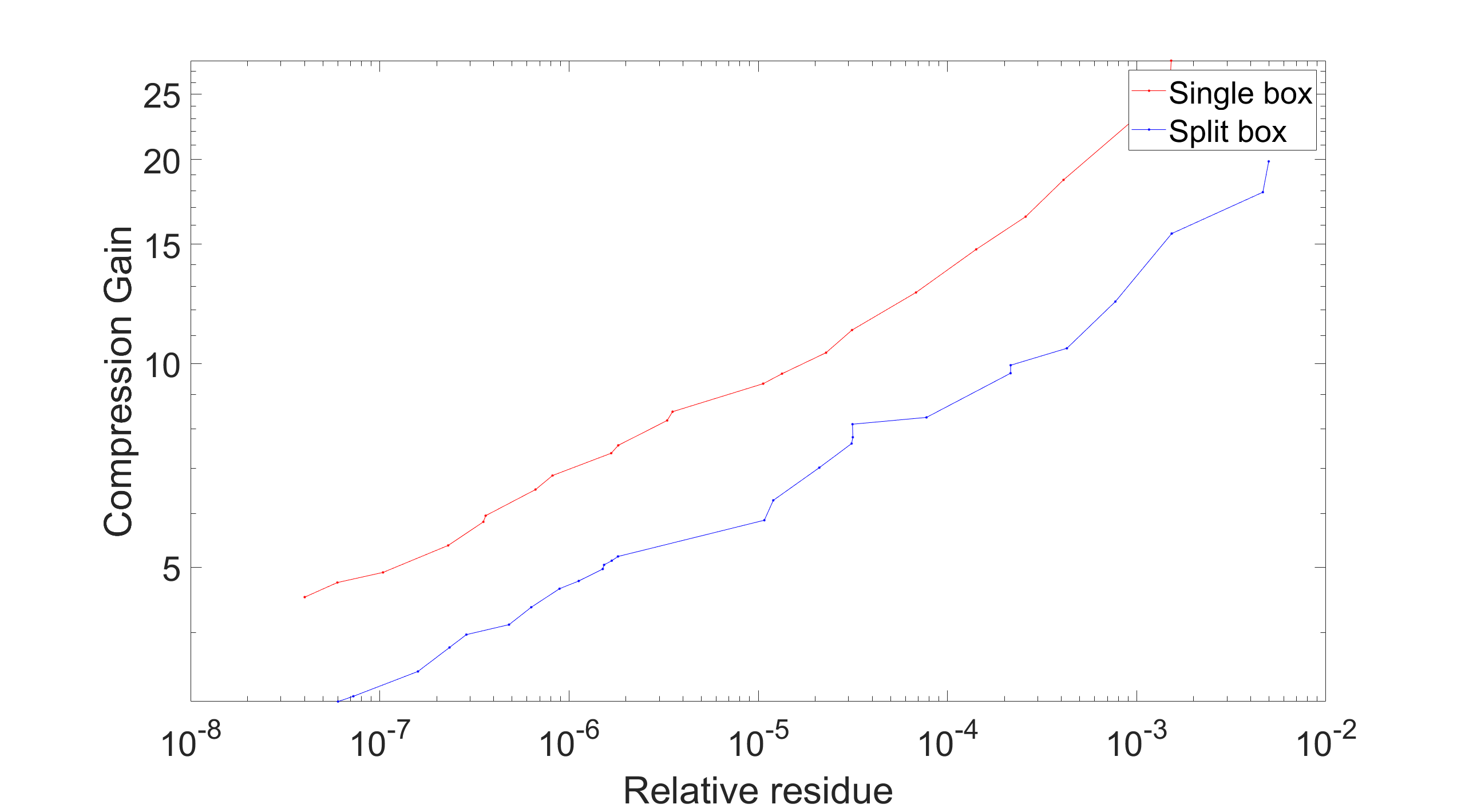}
    \caption{Compression Gain vs the approximation error.}
\label{fig:compgainvsrelres}
\end{figure}

In Figure \ref{fig:compgainvsrelres}, the compression gain defined as in (\ref{eqn_cgain}) is shown for the two different approaches as a function of the relative error obtained during the compression.
It is worth noting that for a given compression gain, the accuracy achieved without splitting the meta-atom is always higher that that with the splitting.

\section{Algorithm for distributing block-block interactions among processors }
\label{appendix:Load}

In the following, the parallelization algorithm is reported. The parallelization is performed over QR-recursive block-block interactions, i.e., meta-atom meta-atom interactions. These interactions are assigned statically to MPI processes using a simple greedy least-loaded rule, and the same mapping is reused both during the assembly phase and for each matrix-vector product.

\begin{tcolorbox}[colframe=black!50, colback=white!10, sharp corners]
    \textbf{Input:}
    
    \hspace{1em} $\mathbf{W_i}$: $i$-th far interaction,  $W_i = (m_i + n_i) r_i$
    
    \hspace{1em} $\mathbf{N_p}$: number of  processes
\end{tcolorbox}

\begin{tcolorbox}[colframe=black!50, colback=white!10, sharp corners]
    \textbf{Output:}  

\hspace{1em} $\mathbf{Wei2Proc}(i)$: $i$-th interaction process owner  

\hspace{1em} $\mathbf{C_k}$: $k$-th process overall load  
\end{tcolorbox}

\begin{algorithmic}

\State $W_i^s \gets \text{sort}(W_i, \text{descending})$

\LComment{Zeroing processor loads}

\For{$s = 1, \dots, N_p$}
    \State $C_s \gets 0$
\EndFor

\LComment{Assign interactions to  processes}
     
    \For{$i = 1, \dots, N_{\text{far}}$}
  
       \State   $k \gets \arg\min_s C_s$          \Comment{Find  least loaded process}    
       
        \State $C_k \gets C_k + W_i^s$          \Comment{Update $k-th$  process load}

         \State $\text{Wei2Proc}(i) \gets k$ \Comment{Assign $W_i^s$  to $k-th$  process}    
    \EndFor

\end{algorithmic}

Differently from existing parallelization algorithm (see \cite{d_mem}), the adopted strategy takes advantage of the specific application allowing for a simple strategy able to achieve good load balance (as demonstrated in Section \ref{sec_eff_p}), while minimization the inter-process communication cost.

\clearpage

\bibliographystyle
{iopart-num}
\bibliography{referencesCF,bibliography}
\end{document}